\documentclass[journal]{IEEEtran}

\newcommand{\beq}{\begin{equation}}
\newcommand{\eeq}{\end{equation}}
\newcommand{\bsq}{\begin{subequations}}
        \newcommand{\esq}{\end{subequations}}
\newcommand{\bq}{\begin{eqnarray}}
\newcommand{\eq}{\end{eqnarray}}
\newcommand{\bqn}{\begin{eqnarray*}}
        \newcommand{\eqn}{\end{eqnarray*}}
\usepackage{bbm}
\usepackage[pdftex]{graphicx}
\usepackage{enumerate}
\usepackage{enumitem} 
	
\usepackage{multirow}
\usepackage{mathptmx}       % selects Times Roman as basic font
\DeclareMathAlphabet{\mathcal}{OMS}{cmsy}{m}{n}
\usepackage{helvet}         % selects Helvetica as sans-serif font
\usepackage{courier}        % selects Courier as typewriter font
\usepackage{type1cm}        % activate if the above 3 fonts are
\usepackage{makeidx}         % allows index generation
\usepackage{graphicx}        % standard LaTeX graphics tool
\usepackage{url} 
\usepackage{multicol}        % used for the two-column index
\usepackage[bottom]{footmisc}% places footnotes at page bottom
\usepackage{ifthen}
\usepackage{subfigure}
\usepackage[usenames,dvipsnames]{color}
\makeindex             % used for the subject index
\usepackage{graphics} % for pdf, bitmapped graphics files
\usepackage{graphicx} % for pdf, bitmapped graphics files
\usepackage{epsfig} % for postscript graphics files
\usepackage{epstopdf}
\usepackage{mathptmx} % assumes new font selection scheme installed
\usepackage{times} % assumes new font selection scheme installed
\usepackage[cmex10]{amsmath}

\usepackage{mathtools}  % loads amsmath
\usepackage{amssymb}  % assumes amsmath package installed
\usepackage{amsthm}

\usepackage{amsfonts}
\usepackage{cite}
\usepackage{verbatim}
\usepackage{comment}
\usepackage{multirow}
\usepackage{cite}
\usepackage{stfloats}
\usepackage{supertabular}
\usepackage{longtable}
\usepackage{tablefootnote}

\DeclareGraphicsExtensions{.pdf,.png,.jpg,.eps,.ps}
\renewcommand{\arraystretch}{1.2}
\usepackage{arydshln}
\usepackage{multicol}
\usepackage{colortbl}
\theoremstyle{definition}

\newtheorem{proposition}{Proposition}

\theoremstyle{definition}
\newtheorem{definition}{Definition}

\ifCLASSINFOpdf
\else
\fi

\usepackage{comment}
\usepackage{ifthen}
\newboolean{showcomments}
\setboolean{showcomments}{true}
\newcommand{\ychen}[1]{\ifthenelse{\boolean{showcomments}}
        { \textcolor{red}{YC: #1}}}
\newcommand{\tongxin}[1]{\ifthenelse{\boolean{showcomments}}
        { \textcolor{blue}{(#1)}}{}}

\usepackage{amsmath}
\allowdisplaybreaks[4]
\usepackage{algorithmicx}

\usepackage[ruled]{algorithm2e}

\begin{document}
        
%       \title{Flexibility Retained Energy Sharing Mechanism among Massive Users}
\title{Improved Approximation of Dispatchable Region in Radial Distribution Networks via Dual SOCP}
        
\author{Yue Chen and
        Changhong Zhao% <-this % stops a space
\thanks{This work was supported by CUHK research startup fund and the Hong Kong Research Grants Council through ECS Award No. 24210220. (Corresponding to C. Zhao)}
\thanks{Y. Chen is with the Department of Mechanical and Automation Engineering, the Chinese University of Hong Kong, HKSAR, China. (email: yuechen@mae.cuhk.edu.hk)}
\thanks{C. Zhao is with the Department of Information Engineering, the Chinese University of Hong Kong, HKSAR, China. (email: chzhao@ie.cuhk.edu.hk)}
}


        \markboth{Journal of \LaTeX\ Class Files,~Vol.~XX, No.~X, Feb.~2019}%
        {Shell \MakeLowercase{\textit{et al.}}: Bare Demo of IEEEtran.cls for IEEE Journals}



        \maketitle
        
        \begin{abstract}
The concept of dispatchable region is useful in quantifying how much renewable generation power the system can handle. In this paper, we aim to provide an improved dispatchable region approximation method in distribution networks. First, based on the nonlinear Dist-Flow model, an optimization problem that minimizes the sum of slack variables is formulated to describe the dispatchable region. The nonconvexity caused by alternating-current (AC) power flow constraints makes it intractable. To deal with this issue, the problem is relaxed to a second-order cone program (SOCP) whose strong dual problem is derived. Then, an SOCP-based projection algorithm is developed to construct a convex polytopic approximation. We prove that the proposed algorithm can generate the accurate SOCP-relaxed dispatchable region under certain conditions. Furthermore, a heuristic method is proposed to approximately remove the regions that make the SOCP relaxation inexact. The final region obtained is the difference of several convex sets and can be nonconvex. Thus, the proposed approach may provide a better approximation of the actually nonconvex dispatchable region than previous work that could construct convex sets only. Numerical results demonstrate that the proposed method can achieve a high accuracy of approximation with simple computation.
        \end{abstract}
        

        % \fi
        \IEEEpeerreviewmaketitle

\section*{Nomenclature}
\addcontentsline{toc}{section}{Nomenclature}
\subsection{Constant parameters}
\begin{IEEEdescription}[\IEEEusemathlabelsep\IEEEsetlabelwidth{${\underline P _{mn}}$,${\overline P _{mn}}$}]
\item[$r_{ij}, \chi_{ij}$]   Resistance and reactance of line $i\rightarrow j$.
\item[$\underline p_i, \overline p_i$]  Controllable active power limit at node $i$.
\item[$\underline q_i, \overline q_i$]  Controllable reactive power limit at node $i$.
\item[$\underline v_i, \overline v_i$]  Voltage safety limit at node $i$.
\item[$\overline \ell_{ij}$]   Current safety limit on line $i\rightarrow j$.
\item[]
\item[$A_f, B_f, A_s$] Constant matrices in the feasibility problem.
\item[$\gamma_f, \gamma_s$] Constant vectors in the feasibility problem.
\item[$A_{y},\!b_{y},\!c_{q},\!\gamma_{q}$] Constant matrices and vectors in SOCP.
\item [$\underline w, \overline w$] Bounds for the initial polytope in Algorithm \ref{alg:approximate-socp}. 
\item [$\delta,\eta, \eta'$] Positive parameters used in Algorithm \ref{alg:remove-inexact}. 
\end{IEEEdescription}

\subsection{Variables}
\begin{IEEEdescription}[\IEEEusemathlabelsep\IEEEsetlabelwidth{${\underline P _{mn}}$,${\overline P _{mn}}$}]
\item[$p_i, q_i$] Controllable power injection at node $i$.
\item[$w_i$] Renewable active power generation at node $i$.
\item[$v_i$] Squared voltage magnitude at node $i$.
\item[$\ell_{ij}$] Squared current magnitude on line $i \rightarrow j$.
\item[$P_{ij}, Q_{ij}$] Active and reactive power flows onto line $i \rightarrow j$.
\item[]
\item[$x$] Vector of state variables $x:=(p,q,v,\ell,P,Q)$.
\item[$z_s$, $z_q$, $\tilde z_q$] Vectors of nonnegative slack variables.
\item[$y_{ij} \in \mathbb{R}^3$] Auxiliary variables in SOCP for line $i\rightarrow j$. 
\item[$\mu_f, \mu_y$] Dual variables for equality constraints. 
\item[$\lambda_s, \lambda_q$] Dual variables for inequality constraints.
\end{IEEEdescription}

\subsection{Optimization problems, values, sets}
\begin{IEEEdescription}[\IEEEusemathlabelsep\IEEEsetlabelwidth{${\underline P _{mn}}$,${\overline P _{mn}}$}]
\item[$\mathrm{FP}(w)$] Feasibility problem for renewable generation $w$.
\item[$\mathrm{fp}(w)$] Minimum objective value of $\mathrm{FP}(w)$.
\item[$\mathcal{W}$] Dispatchable region of $w$, in which $\mathrm{fp}(w)=0$.
\item[$\mathrm{FP}'(w)$] SOCP relaxation of $\mathrm{FP}(w)$.
\item[$\mathrm{fp}'(w)$] Minimum objective value of $\mathrm{FP}'(w)$.
\item[$\mathcal{W}'$] SOCP-relaxed dispatchable region of $w$.
\item[]
\item[$\mathrm{DP}'(w)$] Dual problem of  $\mathrm{FP}'(w)$, also an SOCP.
\item [$D_w(\mu,\lambda)$] Dual objective function. 
\item [$\mathrm{dp}'(w)$] Maximum objective value of $\mathrm{DP}'(w)$.
\item[$\mathrm{DP}''(w,\delta)$] Dual SOCP with feasible set tightened by $\delta$.
\item[$\mathrm{dp}''(w,\delta)$] Maximum objective value of $\mathrm{DP}''(w,
\delta)$.
\item[]
\item[$\mathcal{W}'_{poly}$] Polytopic approximation of $\mathcal{W}'$, by Algorithm \ref{alg:approximate-socp}.
\item[$\mathcal{\tilde W}$] SOCP-inexact region of $w$.
\item[$\mathcal{\tilde W}_d$] Approximation of $\mathcal{\tilde W}$ using dual SOCP.
\item[$\mathcal{\tilde W}_{poly}$] Polytopic approximation of $\mathcal{\tilde W}$, by Algorithm \ref{alg:remove-inexact}.
\end{IEEEdescription}

\section{Introduction}
\IEEEPARstart{W}{ith} the benefits of near-zero carbon emissions and low operating costs, distributed renewable power is experiencing tremendous expansion in recent years \cite{ahmad2018dynamic}. Meanwhile, its volatile and intermittent features pose great challenges to electric grid operation, especially the distribution system. Unlike the bulk system, distribution system has few controllable units and a stronger coupling of active and reactive power flows due to its high resistance to reactance ratio \cite{shen2022admissible}. This makes it even harder to accommodate the fluctuating renewable generations. Therefore, characterizing the renewable power capacities that can be safely hosted by a distribution network prior to its actual operation is vital. This necessitates finding all renewable power outputs that can ensure \emph{solvability} of the power flow equations and satisfaction of \emph{safety limits}.

The first requirement is \emph{solvability} of the power flow equations. For a transmission network modeled by direct-current (DC) power flow, solvability is easy to check since a closed-form solution can be obtained \cite{soroudi2015stochastic}. However, for distribution networks, the lossless DC model is not accurate enough since the distribution lines have higher resistance to reactance ratios. Some literature proved sufficient conditions under which the alternating current (AC) power flow equations are solvable, by utilizing Banach fixed-point theorem for contraction mappings \cite{bolognani2015existence, wang2016explicit} or Brouwer fixed-point theorem for continuous mappings over compact convex sets \cite{dvijotham2017solvability, simpson2017theory}. Nonetheless, those methods cannot be readily applied to output the dispatchable region since they are based on power flow equations and can hardly deal with inequality safety constraints. %Those methods cannot be readily applied to the analysis of dispatchable region, whose model consists not only of equations but inequality safety constraints.

To further take into account the second requirement, i.e., satisfaction of safety limits, optimization based methods were developed. Two well-known concepts are the do-not-exceed limit (DNEL) \cite{zhao2014variable} and the dispatchable region \cite{wei2014dispatchable}. The DNEL provides an allowable power interval for each renewable generator based on robust optimization. Data-driven approach \cite{qiu2016data} and topology control \cite{korad2015enhancement} were incorporated to improve the accuracy of DNEL. The correlation between different renewable generators is ignored in DNEL, so the obtained capacity regions can be conservative. The dispatchable region further considers those correlations and provides the exact region consisting of all renewable power outputs that can be accommodated. An adaptive constraint generation algorithm was proposed to generate the dispatchable region \cite{wei2015real}. The interaction between different prosumers with renewable generators was considered in \cite{chen2021energy}. Similarly, dispatchable region can be applied to quantify the allowable variation of loads based on Fourier–Motzkin elimination \cite{abiri2016loadability}. The above studies are based on DC power flow models.

As mentioned above, AC power flow model is a must for a distribution network. Reference \cite{chen2018convex} solved nonlinear programs to get a set of boundary points that each make a different safety limit binding, and then built a dispatchable region heuristically as the convex hull of those boundary points. Linearized models were used to approximate the real dispatchable region under AC power flow \cite{wan2016maximum, liu2019real}. However, there is no guarantee that all scenarios inside the obtained region are feasible. Reference \cite{shen2022admissible} used the intersection of the dispatchable regions generated from two linearized models to output a more accurate approximation. To guarantee feasibility, certified inner approximations of dispatchable regions were solved from convex programs based on a tightened-relaxed second-order cone approximation \cite{nick2017exact} or refined linear approximations \cite{nguyen2018constructing, nazir2019convex} to AC power flow. 
However, such estimation typically only works for a specific objective function that merely explores the dispatchable region towards a single direction or with a specific shape of the renewable generation vector. Moreover, all the aforementioned regions are convex, while the actual dispatchable region can be nonconvex due to the AC power flow constraints. %This may limit its application scenarios and the information that can be provided.

In this paper, we propose an alternative method to complement the literature above. Our main contributions are two-fold:

\begin{itemize}
\item[1)] \emph{Accurate Dispatchable Region of the Second-Order Cone (SOC) Relaxed Model.} A nonlinear Dist-Flow model based optimization problem is developed to characterize the dispatchable region, which is hard to solve due to its nonconvexity. Therefore, we first relax the problem to a convex second-order cone program (SOCP). Then, unlike reference \cite{shen2022admissible} that further linearized the SOC constraint using polyhedral approximation \cite{ben2001polyhedral}, we generate the dispatchable region directly without further linearization. To be specific, the dual problem of the SOCP is derived and strong duality holds as proven in Proposition \ref{prop:strong-duality}. We then propose a projection algorithm (Algorithm \ref{alg:approximate-socp}) to construct a polytopic approximation of the SOCP-relaxed dispatchable region. We prove that the approximation is accurate under certain conditions.

\item[2)] \emph{Removal of SOCP-Inexact Regions.} The other inaccuracy lies in the possible inexactness of the SOC relaxation. In fact, the actual dispatchable region may be nonconvex, but the algorithms developed in previous studies can only generate convex regions. Distinctly, we propose a heuristic method to find out the SOCP-inexact regions by requiring the corresponding dual variables to be larger than some small positive values. Removing the SOCP-inexact regions from the region generated by Algorithm 1 from the SOC relaxed model, we can build a tighter approximation of the actual dispatchable region. The proposed method provides an innovative idea for constructing an accurate dispatchable region as the difference of several convex sets. Numerical results show that the proposed method can approximate the complicated dispatchable region with a simple polytope after moderate computation, while preserving relatively good accuracy. It can also reach a satisfactory balance between ensuring safety and reducing conservatism.
\end{itemize}

The rest of this paper is organized as follows. Section \ref{sec:model} introduces the power network model we use. Section \ref{sec:problem} defines the dispatchable region and the optimization problem to characterize it. Section \ref{sec:methods} elaborates our method to approximate the dispatchable region. Section \ref{sec:numerical} reports numerical experiments, and Section \ref{sec:conclusion} concludes the paper.   

\section{Power network model}\label{sec:model}

Consider the single-phase equivalent model of a distribution network, which is a radial graph with a set $\mathcal{N}$ of nodes and a set $\mathcal{L}$ of lines. Index the nodes as $\mathcal{N} = \{0, 1,\dots,N\}$, where $0$ represents the root node (slack bus). 
For convenience, we treat the lines as directed; for example, if a line connects nodes $i,j \in \mathcal{N}$, where node $i$ is closer to the root than node $j$, then the line directs from $i$ to $j$ and is denoted by $i \rightarrow j$. The power flow in the network at a particular time instant can be modeled by the classic Dist-Flow equations purely in real numbers \cite{baran1989optimal, farivar2013branch}, elaborated as follows. 

At each node $i\in \mathcal{N}$: let $v_i$ denote the squared voltage magnitude; aggregate all the \emph{controllable} power sources and loads into a complex power injection $p_i + \mathrm{j} q_i$; denote the \emph{uncontrollable} active power generation of a renewable energy source as $w_i$. 
Let $\ell_{ij}$ denote the squared current magnitude through each line $i\rightarrow j$. 
Let $P_{ij}$ and $Q_{ij}$ denote the \textit{net} active and \textit{net} reactive power, respectively, that are sent by node $i$ onto line $i\rightarrow j$; they are different from the net power arriving at node $j$ due to power loss, and are negative if node $i$ receives power from line $i\rightarrow j$. 
Let $r_{ij}$, $\chi_{ij}$ denote the constant resistance and reactance of line $i\rightarrow j$, respectively.
The Dist-Flow equations are:
\begin{subequations}\label{eq:dist-flow}
\begin{IEEEeqnarray}{rrCl}
\forall i\rightarrow j: 
&P_{ij} -  r_{ij} \ell_{ij} - \sum_{k: j\rightarrow k} P_{jk}+ p_j+ w_j &=& 0  \label{eq:dist-flow:p} \\
&Q_{ij} - \chi_{ij} \ell_{ij} - \sum_{k: j\rightarrow k} Q_{jk}+ q_j &=& 0  \label{eq:dist-flow:q} \\
v_{i} -  v_{j}  & - 2(r_{ij} P_{ij} + \chi_{ij} Q_{ij}) + (r_{ij}^2 + \chi_{ij}^2) \ell_{ij} &=& 0 \label{eq:dist-flow:v} \\
& P_{ij}^2 + Q_{ij}^2 - v_i \ell_{ij} &=& 0. \label{eq:dist-flow:s}
\end{IEEEeqnarray}
\end{subequations}

Suppose renewable energy sources only exist at a subset of nodes $\mathcal{N}_w \subseteq \mathcal{N}\backslash \{0\}$, whose cardinality is $W:=|\mathcal{N}_w|$. For nodes $i\notin \mathcal{N}_w$, set constant $w_i\equiv 0$.  
The variables in Dist-Flow equations \eqref{eq:dist-flow} are grouped as follows:
\begin{itemize}
\item Renewable power generation $w:=(w_i,~i\in \mathcal{N}_w)\in \mathbb{R}^W$, which is treated as input to the system;
%\item Dispatchable power injections $(p,q)$ stacked as a column vector, where $p=(p_1,...,p_N)$ and $q=(q_1,...,q_N)$;
\item State variables $x:=(p,q, v, \ell, P, Q)$, where each of $p$, $q$, $v$, $\ell$, $P$, $Q$ is a column vector indexed by $\{1,...,N\}$.      
\end{itemize}

\textit{Remark:} Without loss of generality, we assume there is only one node, indexed as node $1$, connected to the root node $0$. In this case, the power exchange between the distribution network and the upper grid at node $0$ is $p_0 + \mathrm{j} q_0 = P_{01} +\mathrm{j} Q_{01}$, so that it is just considered as part of $(P,Q)$, not $(p,q)$. 
As customary, assume $v_0$ is a given constant and thus not in state variable $v$.  
The radial network has $N$ lines, where each line $i\rightarrow j$ can be uniquely indexed by its destination node $j$, so that we can index line variables $\ell$, $P$, $Q$ by $\{1,...,N\}$.

%Given $(p,q) \in \mathbb{R}^{2N}$, equation \eqref{eq:dist-flow} is a set of $(4N)$ equations with $(4N)$ real variables $x=(v,\ell, P,Q)$. 
Assume known capacity limits of controllable power:
\begin{subequations}\label{eq:pq_limits}
\begin{IEEEeqnarray}{rCl}
\underline p_i \leq p_i \leq \overline p_i, &\quad& \forall  i =1,...,N \label{eq:pq_limits:p} \\
\underline q_i \leq q_i \leq \overline q_i, && \forall i =1,...,N \label{eq:pq_limits:q} 
\end{IEEEeqnarray}
\end{subequations}
At any node $i$ where there are only fixed (or zero) power injections, the constant limits can be set as $\underline p_i = \overline p_i$ ($=0$) and/or $\underline q_i = \overline q_i$ ($=0$).
In addition, power system operations require the following safety limits to be satisfied:
\begin{subequations}\label{eq:vl_limits}
\begin{IEEEeqnarray}{rCl}
\underline v_i \leq v_i \leq \overline v_i, &\quad& \forall  i =1,...,N \label{eq:vl_limits:v} \\
0\leq \ell_{ij} \leq \overline \ell_{ij}, && \forall i\rightarrow j  \label{eq:vl_limits:ell} 
\end{IEEEeqnarray}
\end{subequations}
where the voltage limits $\underline v_i$, $\overline v_i$ for all nodes $i$ and the current limits $\overline \ell_{ij}$ for all lines $i\rightarrow j$ are given as positive constants.

With the model above, we next define and analyze the dispatchable region of renewable power generation. 

\section{Dispatchable region and relaxation}\label{sec:problem}

In this paper, the \emph{dispatchable region} is the region of renewable power generation $w$, for which there is a feasible dispatch. Its formal definition is provided below.

%Consider the net active and/or reactive power injections at some nodes to be known constant numbers such as zero. Let subvector $d \in \mathbb{R}^{D}$ collect all such known constant elements in $(p,q)\in\mathbb{R}^{2N}$, and let $w \in \mathbb{R}^W$ collect all the other elements, i.e., the unknown variable power injections, with $D+U=2N$. In practice, $U$ can be much smaller than $(2N)$ \cite{wei2014dispatchable, liu2019real, chen2021energy}.

\begin{definition}\label{def:feasibility}
A vector of renewable power generation $w\in\mathbb{R}^W$ has a \textbf{feasible dispatch} if there exists $x=(p,q,v,\ell,P,Q) \in \mathbb{R}^{6N}$ such that $(w,x)$ satisfies power flow equations \eqref{eq:dist-flow}, capacity limits \eqref{eq:pq_limits}, and safety limits \eqref{eq:vl_limits}. The \textbf{dispatchable region} of renewable power generation is defined as:
\begin{IEEEeqnarray}{rCl}\nonumber
\mathcal{W}&:=& \left\{w \in \mathbb{R}^W~|~w~\text{has a feasible dispatch.}\right\}
\end{IEEEeqnarray}
\end{definition}

%\textit{Remark:} At every node $i$ with unknown variable power supply $u_i^s:=(p_i^s, q_i^s)$ and/or demand $u_i^d := (p_i^d, q_i^d)$, its net power injection is $u_i = u_i^s - u_i^d$. From the dispatchable region $\mathcal{W}$ of $u$, one can derive the dispatchable region $\mathcal{W}^s$ for supply $u^s$ if the operating region $\mathcal{W}^d$ for demand $u^d$ is known, and vice versa. Actually $\mathcal{W}^s = \mathcal{W} + \mathcal{W}^d$, the Minkowski sum. Oftentimes $\mathcal{W}^d$ is specified as $\underline u^d \leq u^d \leq \overline u^d$ (element-wise); if $\mathcal{W}$ is approximated by a finite union of convex polytopes, as discussed in Section \ref{sec:methods}, then the said Minkowski sum can be efficiently and accurately computed \cite{nazir2018inner}.      

For conciseness, we rewrite the linear part \eqref{eq:dist-flow:p}--\eqref{eq:dist-flow:v} of Dist-Flow equations as
$A_f x \!+\! B_{f} w \!+\! \gamma_f = 0$ and affine inequalities \eqref{eq:pq_limits}--\eqref{eq:vl_limits} as $A_s x + \gamma_s  \leq 0$, where both equality and inequality are element-wise, and constant matrices and vectors $A_f$, $B_f$, $\gamma_f$, $A_s$, $\gamma_s$ are provided in Appendix-A.
Given any $w$, we introduce the following optimization to check its feasibility.%\textit{feasibility problem} as the following optimization: 
\begin{subequations}\label{eq:opt-feasibility}
\begin{IEEEeqnarray}{rCl}
\mathrm{FP}(w):~\min &~& 1^\intercal \tilde z \label{eq:opt-feasibility:obj}
\\ 
\text{over} && x=(p,q,v,\ell,P,Q), ~\tilde z=(z_s, z_q, \tilde z_q)\geq 0 \nonumber
\\
\text{s. t.} && A_f x + B_f w + \gamma_f = 0 \label{eq:opt-feasibility:lin-dist-flow}
\\ &&  A_s x +\gamma_s \leq z_s  \label{eq:opt-feasibility:limits}
\\ && P_{ij}^2 + Q_{ij}^2 - v_i \ell_{ij} \leq z_{q,ij}, ~\forall i\rightarrow j \label{eq:opt-feasibility:nonlinear-small}
\\ && v_i \ell_{ij}-(P_{ij}^2 + Q_{ij}^2)\leq \tilde z_{q,ij}, ~\forall i\rightarrow j \label{eq:opt-feasibility:nonlinear-big}
\end{IEEEeqnarray} 
\end{subequations}
where $1^\intercal$ in objective \eqref{eq:opt-feasibility:obj} is a row vector of all ones. Any element of the slack variable $\tilde z$ can increase as needed to satisfy the corresponding inequality constraint, but only $\tilde z=0$ can guarantee feasibility in terms of \eqref{eq:dist-flow}--\eqref{eq:vl_limits}.   
Therefore, denoting the minimum objective value of $\mathrm{FP}(w)$ as $\mathrm{fp}(w)$, the dispatchable region in Definition \ref{def:feasibility} is equivalently:
\begin{IEEEeqnarray}{rCl}\nonumber
\mathcal{W}&=& \left\{w \in \mathbb{R}^W~|~\mathrm{fp}(w)=0 \right\}.
\end{IEEEeqnarray}

Due to the nonconvex quadratic inequality constraint \eqref{eq:opt-feasibility:nonlinear-big}, problem $\mathrm{FP}(w)$ is nonconvex and thus hard to analyze. 
By removing \eqref{eq:opt-feasibility:nonlinear-big} and rewriting \eqref{eq:opt-feasibility:nonlinear-small}, we relax $\mathrm{FP}(w)$ to a convex \emph{second order cone program (SOCP)}:
\begin{subequations}\label{eq:opt-feasibility-relaxed}
\begin{IEEEeqnarray}{rCl}
\mathrm{FP}'(w):~\min &~& 1^\intercal z \label{eq:opt-feasibility-soc:obj}
\\ 
\text{over} && x,~y,~z=(z_s, z_q)\geq 0 \nonumber
\\
\text{s. t.} && \text{\eqref{eq:opt-feasibility:lin-dist-flow}--\eqref{eq:opt-feasibility:limits}} \nonumber \\
&& y = A_{y} x + b_{y} \label{eq:opt-feasibility:soc-substitute}\\
 \| y_{ij} \|_2 &&\leq c_{q,ij} x + \gamma_{q,ij}+z_{q,ij},~  \forall i\rightarrow j \label{eq:opt-feasibility:soc}
\end{IEEEeqnarray} 
\end{subequations}
where $y\in \mathbb{R}^{3N}$, $A_y \in\mathbb{R}^{(3N)\times(6N)}$, and $b_y \in \mathbb{R}^{3N}$ vertically stack $y_{ij}\in\mathbb{R}^3$,  $A_{y,ij} \in \mathbb{R}^{3\times(6N)}$, and $b_{y,ij}\in \mathbb{R}^3$ respectively for all lines $i\rightarrow j$. Row vector $c_{q,ij} \in \mathbb{R}^{1\times(6N)}$ and scalar number $\gamma_{q,ij} \in \mathbb{R}$ are also stacked vertically for all $i\rightarrow j$ as $c_q \in \mathbb{R}^{N\times (6N)}$ and $\gamma_q \in \mathbb{R}^{N}$. The constant matrices and vectors $A_{y}$, $b_{y}$, $c_{q}$, $\gamma_{q}$ are provided in Appendix-B, which make: 
\begin{IEEEeqnarray}{rCl}
A_{y,ij} x + b_{y,ij} &=& [2P_{ij}, ~2Q_{ij}, ~ v_i \!-\! \ell_{ij}]^\intercal,\quad \forall i\rightarrow j \nonumber \\
c_{q,ij} x + \gamma_{q,ij} &=& v_i + \ell_{ij}, \qquad\qquad\qquad\quad \forall i\rightarrow j  \nonumber
\end{IEEEeqnarray}
and thus make \eqref{eq:opt-feasibility:soc-substitute}--\eqref{eq:opt-feasibility:soc} equivalent to \eqref{eq:opt-feasibility:nonlinear-small}.\footnote{Given $x$, the values of $z_q$ in \eqref{eq:opt-feasibility:nonlinear-small} and \eqref{eq:opt-feasibility:soc} are generally not equal, but we do not differentiate notation due to their identical role as slack variables.} 

Problem $\mathrm{FP}'(w)$ facilitates the definition of an \textit{SOCP-relaxed} dispatchable region:
\begin{IEEEeqnarray}{rCl}\nonumber
\mathcal{W}'&:=& \left\{w \in \mathbb{R}^W~|~\mathrm{fp}'(w)=0 \right\}
\end{IEEEeqnarray}
where $\mathrm{fp}'(w)$ is the minimum objective value of $\mathrm{FP}'(w)$. It is obvious that~$\mathcal{W} \subseteq \mathcal{W}'$, i.e., $\mathcal{W}'$ is a relaxation of $\mathcal{W}$.

A common practice to further simplify the dispatchable-region characterization is to outer approximate the second-order cone \eqref{eq:opt-feasibility:soc} with a polytopic cone, which can achieve arbitrary precision by constructing sufficiently many planes tangent to the surface of the second-order cone \cite{ben2001polyhedral,chen2018energy}. Consequently, $\mathrm{FP}'(w)$ is relaxed to a linear program, and then the algorithm in \cite{wei2014dispatchable,wei2015real, chen2021energy} can be employed to get a convex polytopic outer approximation of $\mathcal{W}'$. 
In this work, we propose an alternative method that does not rely on such linearization. Instead, we work directly on the SOCP $\mathrm{FP}'(w)$ and its dual problem to preserve the intrinsic nonlinearity of the AC power flow model and hence the accuracy of our characterization.

\section{Polytopic approximation algorithms}\label{sec:methods}

To offer a closed-form approximation of dispatchable region $\mathcal{W}$, we first develop a convex polytopic approximation of its relaxation $\mathcal{W}'$ via the dual problem of SOCP $\mathrm{FP}'(w)$. 
We then develop a heuristic method to approximately remove the renewable generations that make the SOCP relaxation inexact, resulting in a tighter approximation of $\mathcal{W}$.   

\subsection{Dual SOCP}

Let $\mu:=(\mu_f, \mu_y)$ denote the dual variables for the equality constraints in problem $\mathrm{FP}'(w)$, with $\mu_f \in \mathbb{R}^{3N}$ for \eqref{eq:opt-feasibility:lin-dist-flow} and $\mu_y\in \mathbb{R}^{3N}$ for \eqref{eq:opt-feasibility:soc-substitute} vertically stacking $\mu_{y,ij}\in \mathbb{R}^3,~\forall i\rightarrow j$.  
Let $\lambda:=(\lambda_s, \lambda_q)$ denote the dual variables for the inequality constraints, with $\lambda_s \in \mathbb{R}^{8N}$ for \eqref{eq:opt-feasibility:limits} and $\lambda_q=(\lambda_{q,ij},~\forall i\rightarrow j) \in \mathbb{R}^N$ for \eqref{eq:opt-feasibility:soc}. Then the Lagrangian of $\mathrm{FP}'(w)$ is:
\begin{IEEEeqnarray}{rCl}\label{eq:lagrangian}
L_u &=& 1^\intercal z \ + \  \mu_f^\intercal (A_f x + B_f w + \gamma_f) \nonumber \\
&& + \lambda_s^\intercal (A_s x + \gamma_s - z_s) +  \mu_{y}^\intercal \left(y -A_{y} x - b_{y}\right) \nonumber \\
&&+  \sum_{i\rightarrow j} \lambda_{q,ij} \left(\| y_{ij} \|_2 -  c_{q,ij} x - \gamma_{q,ij} - z_{q,ij}\right) \nonumber \\ 
&=&  z^\intercal (1-\lambda) + \sum_{i\rightarrow j} \left(y_{ij}^\intercal \mu_{y,ij} +\| y_{ij} \|_2\lambda_{q,ij}\right)\nonumber \\
&&  + x^\intercal \left(A_f^\intercal \mu_f + A_s^\intercal \lambda_s  - A_{y}^\intercal \mu_{y} - c_{q}^\intercal \lambda_{q}\right)   \nonumber \\
&&  + \mu_f^\intercal (B_f w + \gamma_f)+\lambda_s^\intercal \gamma_s  -  \mu_{y}^\intercal b_y -\lambda_{q}^\intercal \gamma_q. \label{eq:Lagrangian}
\end{IEEEeqnarray}

Through $\min_{z\geq 0, x, y} L_u(x,y,z; \mu,\lambda)$ we can get the dual objective function.
By \eqref{eq:Lagrangian}, $L_u$ can only attain a finite minimum over $(z\geq 0, x,y)$ when the dual variables satisfy: 
\begin{subequations}\label{eq:dual-feasibility}
\begin{IEEEeqnarray}{rCl}
0 \leq &~\lambda~ & \leq 1 \label{eq:dual-feasibility:z} \\
A_f^\intercal \mu_f + A_s^\intercal \lambda_s  &=&  A_y^\intercal \mu_y + c_q^\intercal \lambda_q \label{eq:dual-feasibility:x}\\
 \|\mu_{y,ij}\|_2 &\leq& \lambda_{q,ij},\quad\forall i\rightarrow j\label{eq:dual-feasibility:y}
\end{IEEEeqnarray}
\end{subequations}
Note that $\lambda\geq 0$ in \eqref{eq:dual-feasibility:z} is a general requirement for all the dual variables associated with inequality constraints, and \eqref{eq:dual-feasibility:y} must hold by noticing
\begin{IEEEeqnarray}{rCl}
y_{ij}^\intercal \mu_{y,ij} +\| y_{ij} \|_2\lambda_{q,ij} &\geq& \left(\lambda_{q,ij} - \|\mu_{y,ij}\|_2 \right)~ \|y_{ij}\|_2. \nonumber
\end{IEEEeqnarray}
When \eqref{eq:dual-feasibility} is satisfied, all the terms containing $(x,y,z)$ in \eqref{eq:Lagrangian} attain their minimum value zero, and hence we obtain the dual problem for $\mathrm{FP}'(w)$, which is also an SOCP:
\begin{IEEEeqnarray}{rCl}
\mathrm{DP}'(w):~\max_{\mu,\lambda} &\quad& \mu_f^\intercal (B_f w + \gamma_f) + \lambda_s^\intercal \gamma_s -  \mu_{y}^\intercal b_y -\lambda_{q}^\intercal \gamma_q \nonumber
\\
\text{s. t.} && \text{\eqref{eq:dual-feasibility}.} \nonumber
\end{IEEEeqnarray} 
 Let $D_w(\mu,\lambda)$ denote the objective function and $\mathrm{dp}'(w)$ denote the maximum objective value of $\mathrm{DP}'(w)$.
The following result lays the foundation for approximating the SOCP-relaxed dispatchable region $\mathcal{W}'$ via the dual SOCP $\mathrm{DP}'(w)$.

\begin{proposition}\label{prop:strong-duality}
For all $w\in\mathbb{R}^W$, strong duality holds between $\mathrm{FP}'(w)$ and $\mathrm{DP}'(w)$, i.e., their optimal values $\mathrm{fp}'(w) = \mathrm{dp}'(w)$. 
\end{proposition}
\begin{IEEEproof}
Consider an arbitrary $w\in\mathbb{R}^W$. Since problem $\mathrm{FP}'(w)$ is convex, it is sufficient to prove Slater's condition \cite[Section 5.2.3]{boyd2004convex}, i.e., existence of $(z\geq 0,x,y)$ that satisfies affine constraints \eqref{eq:opt-feasibility:lin-dist-flow}\eqref{eq:opt-feasibility:limits}\eqref{eq:opt-feasibility:soc-substitute} and strictly satisfies \eqref{eq:opt-feasibility:soc}. 

Indeed, it is adequate to find a point $x=(p,q,v,\ell, P, Q)$ to satisfy \eqref{eq:opt-feasibility:lin-dist-flow}, i.e., \eqref{eq:dist-flow:p}--\eqref{eq:dist-flow:v}; then one can explicitly determine $y$ by \eqref{eq:opt-feasibility:soc-substitute} and always find large enough $z$ to make \eqref{eq:opt-feasibility:limits}\eqref{eq:opt-feasibility:soc} (strictly) feasible, satisfying Slater's condition. 
Such a point $x$ can be easily found as follows: set $p=q=\ell=0 \in \mathbb{R}^N$; determine $(P,Q)$ backward from the leaves to the root of the radial network, using \eqref{eq:dist-flow:p}--\eqref{eq:dist-flow:q}; then determine $v$ forward from the root to the leaves, using \eqref{eq:dist-flow:v}. This completes the proof.
\end{IEEEproof}

By Proposition \ref{prop:strong-duality}, the relaxed region $\mathcal{W}'$ is equivalently:
\begin{IEEEeqnarray}{rCl}
&&\mathcal{W}'= \left\{w \in \mathbb{R}^W~|~\mathrm{dp}'(w)=0 \right\}  \nonumber\\
&=& \left\{w \in \mathbb{R}^W~|D_w(\mu,\lambda)\leq 0,~\forall (\mu,\lambda)~\text{satisfying \eqref{eq:dual-feasibility}}\right\} \label{eq:relaxed-region}
\end{IEEEeqnarray}
where the second equality holds because $D_w(\mu,\lambda)=0$ can always be attained at the dual feasible point $(\mu,\lambda)=0$.

\begin{proposition}\label{proposition:convexityU}
$\mathcal{W}'$ is a convex set.
\end{proposition}
\begin{IEEEproof}
Consider arbitrary $w_1, w_2 \in \mathcal{W}'$ and $t \in [0,1]$. Denote $w_t:=t w_1 + (1-t) w_2$. Then for every $(\mu,\lambda)$ satisfying \eqref{eq:dual-feasibility}, we have:
\begin{IEEEeqnarray}{rCl}
D_{w_t}(\mu, \lambda) 
&=& t D_{w_1} (\mu,\lambda)  + (1-t)D_{w_2} (\mu,\lambda) \nonumber \\
&\leq& t\cdot 0 + (1-t) \cdot 0 = 0 \nonumber
\end{IEEEeqnarray}
where the first equality is due to linearity of $D_w(\mu,\lambda)$ with respect to $w$ when $(\mu, \lambda)$ is fixed, and the inequality holds because $w_1, w_2 \in \mathcal{W}'$. Therefore $w_t \in \mathcal{W}'$. By the definition of a convex set, $\mathcal{W}'$ is convex.
\end{IEEEproof}

\subsection{Approximating SOCP-relaxed dispatchable region}

\begin{algorithm}[t]\label{alg:approximate-socp}
	\caption{Approximate $\mathcal{W}'$}
	1. \textbf{Initialization:} $\mathcal{W}'_{poly} = \left\{w \in \mathbb{R}^W~|~\underline w \leq w \leq \overline w\right\}$ for sufficiently low $\underline w$ and high $\overline w$; $\mathcal{V}_{safe}=\emptyset$; $c=0$.
	
	2. Update vertex set $vert\left(\mathcal{W}'_{poly}\right)$. Let $\mathrm{dp}'_{max}=0$; 
	
	\For{$w\in \mbox{vert}\left(\mathcal{W}'_{poly}\right)$ and $w\notin \mathcal{V}_{safe}$}{
		solve $\mathrm{DP}'(w)$ to obtain an optimal solution $(\mu^*, \lambda^*)$ and maximum objective value $\mathrm{dp}'(w)$;
		
		\lIf{$\mathrm{dp}'(w)>\mathrm{dp}'_{max}$}{
		   
		   \qquad $\mathrm{dp}'_{max} \leftarrow \mathrm{dp}'(w)$;
		
		   \qquad $(\mu_{max},\lambda_{\max}) \leftarrow (\mu^*,\lambda^*)$
		}
		\lElseIf{$\mathrm{dp}'(w)\leq 0$}{
		   $\mathcal{V}_{safe} = \mathcal{V}_{safe} \cup \{w\}$
		}
	}
	\eIf{$\mathrm{dp}'_{max} = 0$ or $c=C_{max}$}{
	    return $\mathcal{W}'_{poly}$.
	}{
	    add to $\mathcal{W}'_{poly}$ a cutting plane:
	    $\mu_{f,max}^\intercal (B_f w + \gamma_f) + \lambda_{s,max}^\intercal \gamma_s \leq  \mu_{y,max}^\intercal b_y +\lambda_{q,max}^\intercal \gamma_q$;
	    
	    $c \leftarrow c+1$;
	    
	    go back to Line 2;
	}
\end{algorithm}

We propose Algorithm \ref{alg:approximate-socp} to approximate $\mathcal{W}'$ defined in \eqref{eq:relaxed-region}. It starts with a region $\mathcal{W}'_{poly}$ that is large enough to contain $\mathcal{W}'$. Then it solves the dual SOCP $\mathrm{DP}'(w)$ for every vertex $w$ of polytope $\mathcal{W}'_{poly}$, records the vertex that most severely violates the condition in \eqref{eq:relaxed-region}, and adds a corresponding cutting plane to remove that vertex from $\mathcal{W}'_{poly}$. Meanwhile, all the vertices that satisfy the condition in \eqref{eq:relaxed-region} are added to $\mathcal{V}_{safe}$ and never checked again. 

\begin{proposition}\label{prop:polytopic-outer-approximation}
The output $\mathcal{W}'_{poly}$ in an arbitrary iteration of Algorithm \ref{alg:approximate-socp} is an outer approximation of $\mathcal{W}'$.  
\end{proposition}
\begin{IEEEproof}
Note the initial $\mathcal{W}'_{poly}$ contains $\mathcal{W}'$. 
We next prove that any cutting plane added to $\mathcal{W}'_{poly}$ would not remove any point in $\mathcal{W}'$.
To show that, consider an arbitrary $w$ removed by a cutting plan whose coefficients are $(\mu_{max},\lambda_{max})$. Then there must be $D_w(\mu_{max},\lambda_{max}) > 0$. Since $(\mu_{max},\lambda_{max})$ is dual feasible satisfying \eqref{eq:dual-feasibility}, we have $w\notin \mathcal{W}'$ by \eqref{eq:relaxed-region}. 
\end{IEEEproof}

 Unlike \cite{wei2014dispatchable, chen2021energy} that based on linear programs, the SOCP-relaxed dispatchable region $\mathcal{W}'$ may not be the intersection of a finite number of cutting planes (i.e., a convex polytope). 
 Therefore, Algorithm \ref{alg:approximate-socp} may not guarantee $\mathrm{dp}'(w)= 0$ for all vertices $w\in vert\left(\mathcal{W}'_{poly}\right)$ in a finite number of iterations. However, if it does so, as what always happens in our numerical experiments, it will produce a nice result as follows.
 
\begin{proposition}\label{prop:polytopic-approximation}
If Algorithm \ref{alg:approximate-socp} terminates with $\mathrm{dp}'_{max} = 0$ in a finite number of iterations, it returns the accurate SOCP-relaxed dispatchable region, i.e., $\mathcal{W}'_{poly}=\mathcal{W}'$.  
\end{proposition}
\begin{IEEEproof}
Proposition \ref{prop:polytopic-outer-approximation} has shown $\mathcal{W}' \subseteq \mathcal{W}'_{poly}$.
If Algorithm \ref{alg:approximate-socp} terminates with $\mathrm{dp}'_{max} = 0$ after adding a finite number of cutting planes, then it returns a convex polytope $\mathcal{W}'_{poly}$. Moreover, all the vertices $w\in vert\left(\mathcal{W}'_{poly}\right)$ satisfy $\mathrm{dp}'(w)= 0$, therefore, $w\in \mathcal{W}'$ by \eqref{eq:relaxed-region}. 
This fact, together with the convexity of $\mathcal{W}'$ shown in Proposition \ref{proposition:convexityU}, implies $\mathcal{W}'_{poly} \subseteq \mathcal{W}'$. Thus we have proved $\mathcal{W}'_{poly}=\mathcal{W}'$.
\end{IEEEproof}

An immediate corollary of Proposition \ref{prop:polytopic-approximation} is that if $\mathcal{W}'$ is not a polytope, then Algorithm \ref{alg:approximate-socp} cannot terminate in a finite number of iterations with $\mathrm{dp}'_{max} = 0$. 
If that happens, one can terminate Algorithm \ref{alg:approximate-socp} when reaching the maximum number of iterations $C_{max}$, to obtain a convex polytopic outer approximation of $\mathcal{W}'$.
In this sense, the outcome of Algorithm \ref{alg:approximate-socp} serves as a posterior indicator of the structure of $\mathcal{W}'$.

\subsection{Removing SOCP-inexact renewable generations}
     
Remember our goal is to characterize the dispatchable region $\mathcal{W}$, whereas $\mathcal{W}'$ studied so far is just a SOCP-relaxation of $\mathcal{W}$. To overcome this drawback, we design a heuristic to approximately remove the SOCP-inexact region $\mathcal{\tilde W}:=\mathcal{W}' \backslash \mathcal{W}$ from $\mathcal{W}'$.
The renewable generations $w \in \mathcal{\tilde W}$ are feasible in terms of the SOCP relaxation $\mathrm{FP}'(w)$ but infeasible in terms of $\mathrm{FP}(w)$, as formally defined below.

\begin{definition}
A vector of renewable power generation $w\in\mathcal{W}'$ is \textbf{SOCP-inexact}, if every optimal solution of $\mathrm{FP}'(w)$ satisfies: 
\begin{IEEEeqnarray}{rCl}
\| y_{ij} \|_2 &<& c_{q,ij} x + \gamma_{q,ij}\quad\text{for some}~ i\rightarrow j. \nonumber
\end{IEEEeqnarray}
The \textbf{SOCP-inexact region} of $w$ is defined as:
\begin{IEEEeqnarray}{rCl}
 \mathcal{\tilde W} &=& \left\{w \in \mathcal{W}'~|~w~\text{is SOCP-inexact}\right\}. \nonumber
\end{IEEEeqnarray}
\end{definition}

Our next focus is to build an approximation of $\mathcal{\tilde W}$. 
For that, we consider the following set defined on the dual SOCP:
\begin{IEEEeqnarray}{rCl}
\mathcal{\tilde W}_d &:=& \{w \in \mathcal{W}'~|~\text{Every optimal solution of}~\mathrm{DP}'(w) \nonumber \\
&&\qquad\qquad\quad \text{satisfies}~\lambda_{q,ij}=0~\text{for some}~i\rightarrow j \}. \nonumber
\end{IEEEeqnarray}
By complementary slackness \cite[Section 5.5.2]{boyd2004convex}, for every primal-dual optimal of $\mathrm{FP}'(w)$ and $\mathrm{DP}'(w)$, there is:
\begin{IEEEeqnarray}{rCl}
\lambda_{q,ij}\left(\| y_{ij} \|_2 - c_{q,ij} x - \gamma_{q,ij}\right)&=&0, \quad \forall i\rightarrow j. \nonumber
\end{IEEEeqnarray}
This implies $\mathcal{\tilde W}\subseteq \mathcal{\tilde W}_d$. Although $\mathcal{\tilde W}= \mathcal{\tilde W}_d$ may not hold, their difference can only occur under rare circumstances where $\lambda_{q,ij} = \| y_{ij} \|_2 \!-\! c_{q,ij} x \!-\! \gamma_{q,ij} = 0$ at a primal-dual optimal. Hence we focus on $\mathcal{\tilde W}_d$ as an approximation of $\mathcal{\tilde W}$. 

Given an arbitrary $w \in \mathcal{\tilde W}_d \subseteq \mathcal{W}'$, the maximum objective value of $\mathrm{DP}'(w)$ is $\mathrm{dp}'(w)=0$ but with some $\lambda_{q,ij}=0$ so the SOC relaxation is inexact (except for some very rare case). To approximate $\tilde{\mathcal{W}}_d$, first we add the following constraint to tighten the dual feasible set \eqref{eq:dual-feasibility}:
\begin{IEEEeqnarray}{rCl}\label{eq:add-dual-constraint}
\lambda_{q} &\geq& \delta
\end{IEEEeqnarray}
where the inequality is element-wise and $\delta\in\mathbb{R}_+^{9N}$ is a vector of all \textit{strictly positive} parameters, whose design will be elaborated later. Consider the tightened dual SOCP:
\begin{IEEEeqnarray}{rCl}
\mathrm{DP}''(w,\delta):~\max_{\mu,\lambda} &\quad&  \mu_f^\intercal (B_f w + \gamma_f) + \lambda_s^\intercal \gamma_s - \mu_{y}^\intercal b_y - \lambda_{q}^\intercal \gamma_q \nonumber
\\
\text{s. t.} && \text{\eqref{eq:dual-feasibility}, \eqref{eq:add-dual-constraint}} \nonumber
\end{IEEEeqnarray} 
and let $\mathrm{dp}''(w,\delta)$ denote its maximum objective value. For $w\in\mathcal{\tilde W}_d$, there must be $\mathrm{dp}''(w,\delta) < 0$, because otherwise $\mathrm{DP}'(w)$ would have an optimal solution that satisfies \eqref{eq:add-dual-constraint}, contradicting the definition of $\mathcal{\tilde W}_d$. 
Actually $\mathrm{dp}''(w,\delta) \leq -\eta$ for some $\eta>0$ that depends on $w$ and $\delta$.

\begin{algorithm}[t]\label{alg:remove-inexact}
	\caption{Approximate $\mathcal{\tilde W}_d$ (or SOCP-inexact $\mathcal{\tilde W}$)}
	1. \textbf{Initialization:} $\mathcal{\tilde W}_{poly} = \mathcal{W}'_{poly}$ returned by Alg. \ref{alg:approximate-socp}. Given positive $\delta$, $\eta$, $\eta'$; $\mathcal{V}_{safe}=\emptyset$; $c=0$;
	
	2. Update vertex set $vert\left(\mathcal{\tilde W}_{poly}\right)$. Let $\mathrm{dp}''_{max}= -\eta$;
	
	\For{$w\in vert\left(\mathcal{\tilde W}_{poly}\right)$ and $w\notin \mathcal{V}_{safe}$}{
		solve $\mathrm{DP}''(w,\delta)$ to obtain an optimal solution $(\mu^*, \lambda^*)$ and maximum objective value $\mathrm{dp}''(w,\delta)$;
		
		\lIf{$\mathrm{dp}''(w,\delta)>\mathrm{dp}''_{max}$}{
		   
		   \qquad $\mathrm{dp}''_{max} \leftarrow \mathrm{dp}''(w,\delta)$;
		
		   \qquad $(\mu_{max},\lambda_{\max}) \leftarrow (\mu^*,\lambda^*)$
		}
		\lElseIf{$\mathrm{dp}''(w,\delta)\leq -\eta$}{
		   $\mathcal{V}_{safe} = \mathcal{V}_{safe} \!\cup\! \{w\}$
		}
	}
	\eIf{$\mathrm{dp}''_{max} = -\eta$ or $c=C_{max}$}{
	    return $\mathcal{\tilde W}_{poly}$.
	}{
	    add to $\mathcal{\tilde W}_{poly}$ a cutting plane:
	    $\mu_{f,max}^\intercal  (B_f w +  \gamma_f) + \lambda_{s,max}^\intercal \gamma_s \leq \mu_{y,max}^\intercal b_y +\lambda_{q,max}^\intercal \gamma_q  -  \eta'$; 
	    
	    $c\leftarrow c+1$;
	   
	    go back to Line 2;
	}
\end{algorithm}

The idea above inspires us to approximate $\mathcal{\tilde W}_d$ (or $\mathcal{\tilde W}$) by
\begin{IEEEeqnarray}{rCl}
\tilde{\mathcal{W}}_d \approx \left\{w \in \mathbb{R}^W|D_w(\mu,\lambda)\leq -\eta,\forall (\mu,\lambda)~\text{satisfying \eqref{eq:dual-feasibility},\eqref{eq:add-dual-constraint}}\right\}  \nonumber
\end{IEEEeqnarray}
To this end, Algorithm \ref{alg:remove-inexact} can be designed using a similar procedure to Algorithm \ref{alg:approximate-socp}.
%The idea above inspires our design of Algorithm \ref{alg:remove-inexact} to approximate $\mathcal{\tilde W}_d$ (or $\mathcal{\tilde W}$), which follows a similar procedure to Algorithm \ref{alg:approximate-socp}.
Algorithm \ref{alg:remove-inexact} returns a convex polytope $\mathcal{\tilde W}_{poly} \subseteq \mathcal{W}'_{poly}$ that guarantees $\mathrm{dp}''(w,\delta) \leq -\eta < 0$ for all $w \in \mathcal{\tilde W}_{poly}$, which is an approximation of $\tilde{\mathcal{W}}_d$ (or $\tilde{\mathcal{W}}$). Removing $\tilde{\mathcal{W}}_{poly}$ from $\mathcal{W}'_{poly}$, we can obtain an approximation $\mathcal{W}_{poly}=\mathcal{W}'_{poly} \backslash \tilde{\mathcal{W}}_{poly}$ of the actual dispatchable region $\mathcal{W}$. To make Algorithm \ref{alg:remove-inexact} more robust, we may choose $\eta'> \eta$ for the added cutting plane in each iteration.

%If Algorithm \ref{alg:remove-inexact} terminates with $\mathrm{dp}''_{max} = -\eta$ in a finite number of iterations, then it returns a convex polytope $\mathcal{\tilde W}_{poly} \subseteq \mathcal{W}'_{poly}$ that guarantees $\mathrm{dp}''(w,\delta) \leq -\eta < 0$ for all $w \in \mathcal{\tilde W}_{poly}$. To make Algorithm \ref{alg:remove-inexact} more robust, we may choose $\eta'> \eta$ for the added cutting plane in each iteration. 

\textit{Remark:} The parameters $\delta$ and $\eta$ are essential for Algorithm \ref{alg:remove-inexact}. A general guideline is that (1) given $\delta$, choosing a smaller $\eta$ and (2) given $\eta$, choosing a bigger $\delta$ will both make $\mathcal{\tilde W}_{poly}$ bigger and lead to a smaller (more conservative) approximation of $\mathcal{W} = \mathcal{W}' \backslash \mathcal{\tilde{W}}$. Moreover, sometime it is difficult for Algorithm \ref{alg:remove-inexact} to use a single convex polytope $\mathcal{\tilde W}_{poly}$ to accurately approximate the most likely nonconvex $\mathcal{\tilde{W}}$.
To deal with this difficulty, we propose to run Algorithm \ref{alg:remove-inexact} multiple times with different vectors $\delta$. As a result, we obtain multiple convex polytopes whose union serves as a better approximation of $\mathcal{\tilde{W}}$. Those vectors $\delta$ can be selected in the following way. We traverse the vertices of $\mathcal{W}'_{poly}$, select one vertex $w$, and solve the dual SOCP $\mathrm{DP}'(w)$ to get an optimal solution $(\mu^*,\lambda^*)$. Then $\delta$ is constructed by keeping all the strictly positive elements of $\lambda_q^*$ as they are, and add a small positive perturbation to all the zero elements.  

%how fast $\mathrm{dp}''(u,\delta)$ decreases with $\delta$ is unknown and different for each $u\in \mathcal{\tilde W}_d$. Therefore, with a specific $(\delta,\eta)$ adopted, it is hard to tell whether Algorithm \ref{alg:remove-inexact} will return an inner or outer approximation of $\mathcal{\tilde W}_d$, or neither.

\begin{table}[t]
        \renewcommand{\arraystretch}{1.3}
        \renewcommand{\tabcolsep}{1em}
        \centering
        \caption{Summary of different regions}
        \label{tab:summary}
        \begin{tabular}{m{1.2cm}<{\centering}m{1.5cm}<{\centering}m{0.1cm}<{\centering}m{1.2cm}<{\centering}m{0.1cm}<{\centering}m{1.5cm}<{\centering}}
                \hline 
           & SOCP-relaxed Region & = & SOCP-exact Region & + & SOCP-inexact Region \\
           \hline
          Actual &  $\mathcal{W}'$ & = & $\mathcal{W}$ & + & $\tilde{\mathcal{W}} \approx \tilde{\mathcal{W}}_d$ \\
          Approx. & $\mathcal{W}'_{poly}$ & =  & $\mathcal{W}_{poly}$ & + & $\tilde{\mathcal{W}}_{poly}$ \\
          Method & Algorithm \ref{alg:approximate-socp} & & $\mathcal{W}'_{poly} \backslash \tilde{\mathcal{W}}_{poly}$ & &  Algorithm \ref{alg:remove-inexact}\\
                \hline
        \end{tabular}
\end{table}

\textit{Summary}. The relationship of different regions mentioned in this paper is summarized in TABLE \ref{tab:summary}. As discussed, the dispatchable region $\mathcal{W}=\mathcal{W}' \backslash \mathcal{\tilde W}$, where $\mathcal{W}'$ is the SOCP-relaxed dispatchable region and $\mathcal{\tilde W}$ is the SOCP-inexact region. 
We develop Algorithm \ref{alg:approximate-socp} to get $\mathcal{W}'_{poly}$, a convex polytopic approximation of $\mathcal{W}'$; and Algorithm \ref{alg:remove-inexact} to get $\mathcal{\tilde W}_{poly}$, a convex polytopic approximation of $\mathcal{\tilde W}$. Algorithm \ref{alg:remove-inexact} can run multiple times to obtain a more accurate approximation of nonconvex $\mathcal{\tilde W}_d$ (or $\mathcal{\tilde W}$). The outputs of multiple runs of Algorithm \ref{alg:remove-inexact} are then removed from $\mathcal{W}'_{poly}$ to obtain a generally nonconvex polytopic approximation of $\mathcal{W}$.
%It can be expressed as a finite union of convex polytopes, which facilitates the computation of Minkowski sum, as remarked after Definition \ref{def:feasibility}. 

%Note that $\mathcal{W}_{poly}$ may not be convex. Indeed, the more accurately $\mathcal{W}_{poly}$ approximates $\mathcal{W}$, the more likely it is nonconvex, since it is well known that the dispatchable region $\mathcal{W}$ is generally nonconvex under AC models. 
%To echo the remark on Minkowski sum after Definition \ref{def:feasibility}, one can express $\mathcal{W}_{poly}$ as a finite union of convex polytopes in the following way. 
% In every iteration of Algorithm \ref{alg:remove-inexact}, consider the half-space opposite to that in Line 13:
% \begin{IEEEeqnarray}{rCl}
% \mu_{f,max}^\intercal (B_f w \!+\! \gamma_f) +\lambda_{s,max}^\intercal \gamma_s &>& \mu_{y,max}^\intercal b \!+\!\lambda_{q,max}^\intercal \gamma -\eta \nonumber 
% \end{IEEEeqnarray}
% and take the intersection of this half-space with $\mathcal{W}'_{poly}$ to obtain a convex polytope. Then $\mathcal{W}_{poly}$ is the union of such convex polytopes over all the iterations till Algorithm \ref{alg:remove-inexact} terminates.

\section{Case Studies} \label{sec:numerical}
\begin{figure}
	\includegraphics[width=0.35\textwidth]{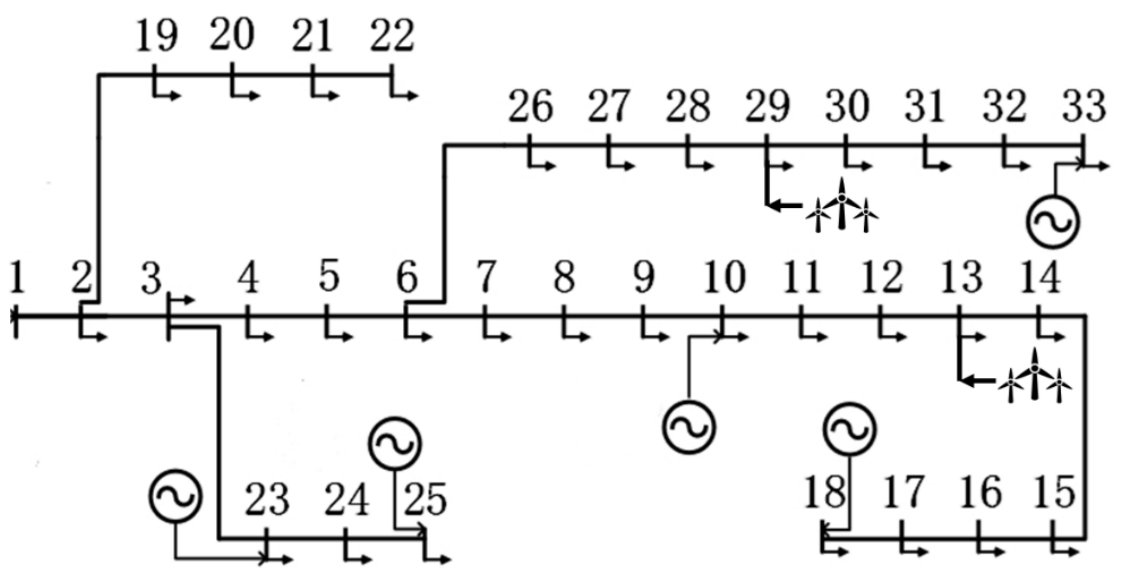}
	\centering
	\caption{IEEE 33-node network model from \cite{chen2018energy}.}
	\label{fig:network}
\end{figure}
In this section, we conduct numerical experiments on the IEEE 33-bus system whose topology is in Fig. \ref{fig:network}. The proposed algorithms are implemented to approximate the dispatchable region of renewable generation $(w_1, w_2)$ at nodes 13 and 29, respectively. Then, we test the impact of several factors and compare with other approaches.
\subsection{Benchmark} \label{sec:numerical-1}
\begin{table}[t]
        \renewcommand{\arraystretch}{1.3}
        \renewcommand{\tabcolsep}{1em}
        \centering
        \caption{Parameters of generators in Benchmark}
        \label{tab:generator}
        \begin{tabular}{cccccc}
                \hline 
                Generator & Location & $\underline{p}_i$ (p.u.) & $\overline{p}_i$ (p.u.)\\
                \hline
                G1 & node 10 & 0.4 & 0.6 \\
                G2 & node 18 & 0.3 & 0.4\\
                G3 & node 23 & 0.4 & 0.6\\
                G4 & node 25 & 0.3 & 0.5\\
                G5 & node 33 & 0.4 & 0.6\\
                \hline
        \end{tabular}
\end{table}
\begin{figure}[t]
        \centering
        \includegraphics[width=0.8\columnwidth]{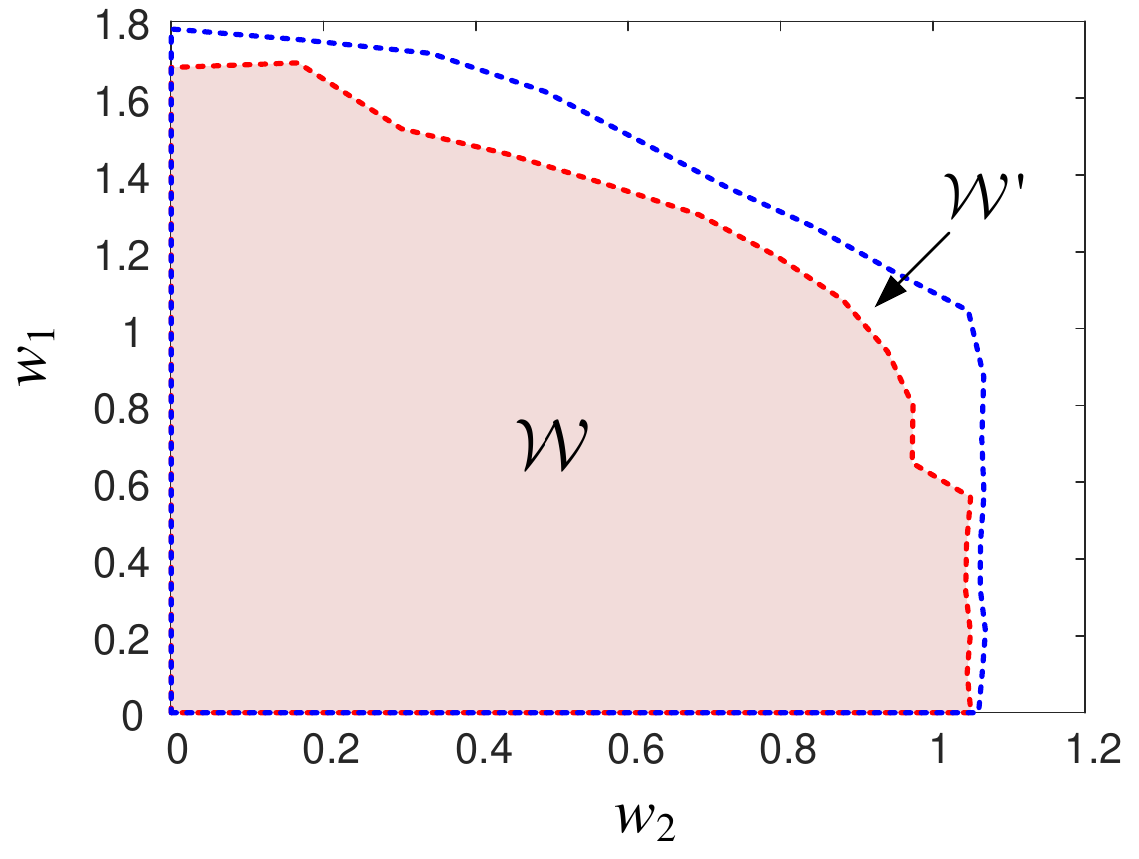}
        \caption{The SOCP-relaxed region $\mathcal{W}'$ (outside) and the actual dispatchable region $\mathcal{W}$ obtained by checking sampled points in the $(w_1,w_2)$ space.}
        \label{fig:actual-relax}
\end{figure}

\begin{figure}[t]
        \centering
        \includegraphics[width=1.0\columnwidth]{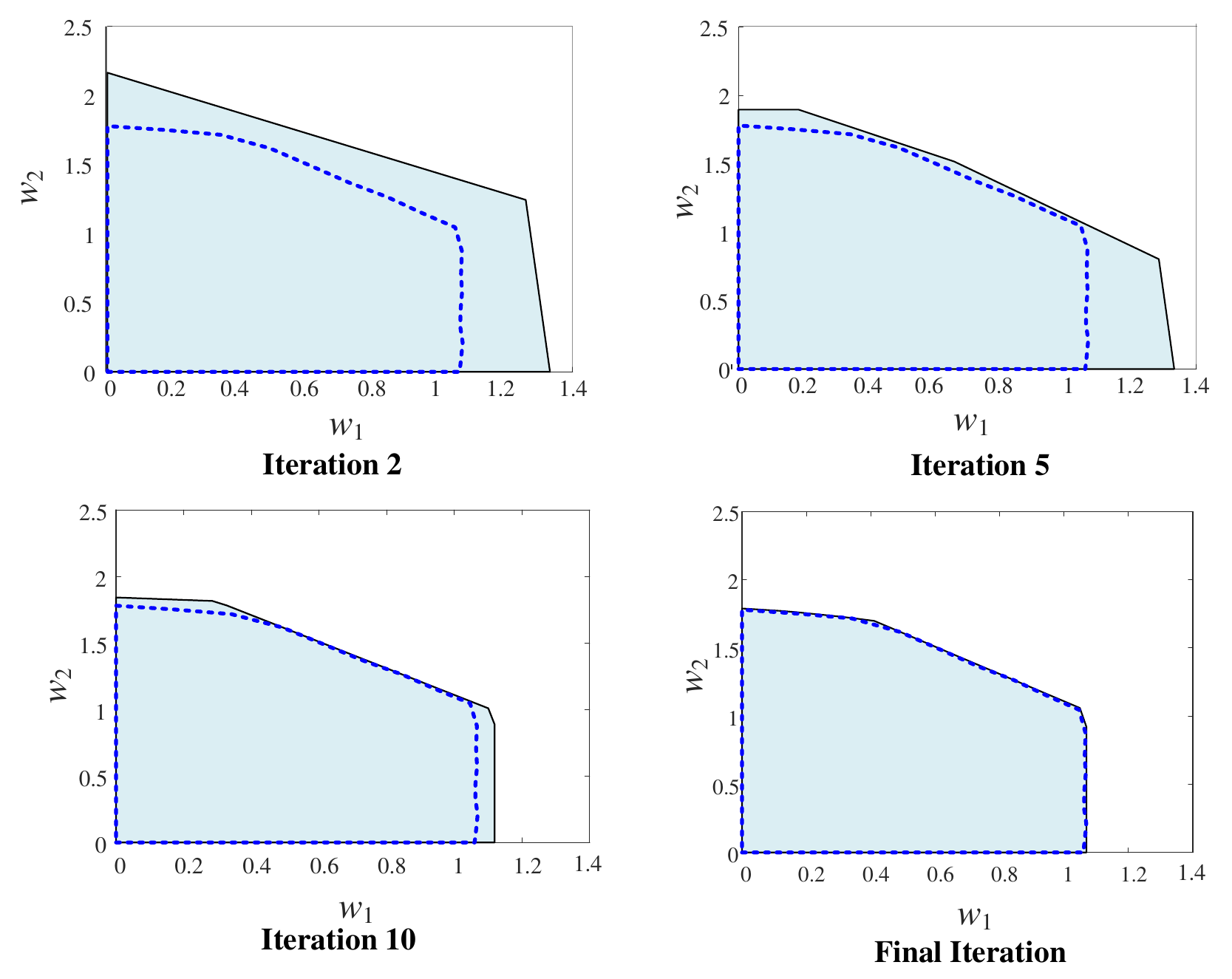}
        \caption{The output of Algorithm \ref{alg:approximate-socp} in different iterations (solid line) to approximate the SOCP-relaxed dispatchable region $\mathcal{W}'$ (dashed line).}
        \label{fig:alg1}
\end{figure}
In the IEEE-33 bus system, there are 5 controllable generators whose parameters are given in TABLE \ref{tab:generator}. Two renewable generators are connected to nodes 13 and 29, respectively.
For comparison, the actual SOCP-relaxed region $\mathcal{W}'$ and the actual dispatchable region without relaxation $\mathcal{W}$ are generated as in Fig. \ref{fig:actual-relax}. This can be done by checking the feasiblity of a nonlinear optimization with \eqref{eq:dist-flow}-\eqref{eq:vl_limits} as its constraints, over sample points $w$ in the $(w_1,w_2)$ space using the nonlinear solver IPOPT. As we can see from Fig. \ref{fig:actual-relax}, the actual dispatchable region $\mathcal{W}$ can be nonconvex and the SOCP-relaxed region is not accurate enough. In the following, we apply the proposed algorithms to output a more accurate region.

First, we test the performance of Algorithm \ref{alg:approximate-socp}. We observe that algorithm terminates with $\mathrm{dp}'_{max} = 0$ in 25 iterations, taking about 289.84s. The output regions $\mathcal{W}'_{poly}$ in the 2nd, 5th, 10th, and final iterations are given in Fig. \ref{fig:alg1}. The Algorithm \ref{alg:approximate-socp} removes the nondispatchable regions iteratively (the blue region is becoming smaller), and finally returns a convex polytope $\mathcal{W}'_{poly}$ exactly the same as the actual SOCP-relaxed region $\mathcal{W}'$ (dashed line). This validates Proposition \ref{prop:polytopic-approximation}.

\begin{figure}[t]
        \centering
        \includegraphics[width=0.8\columnwidth]{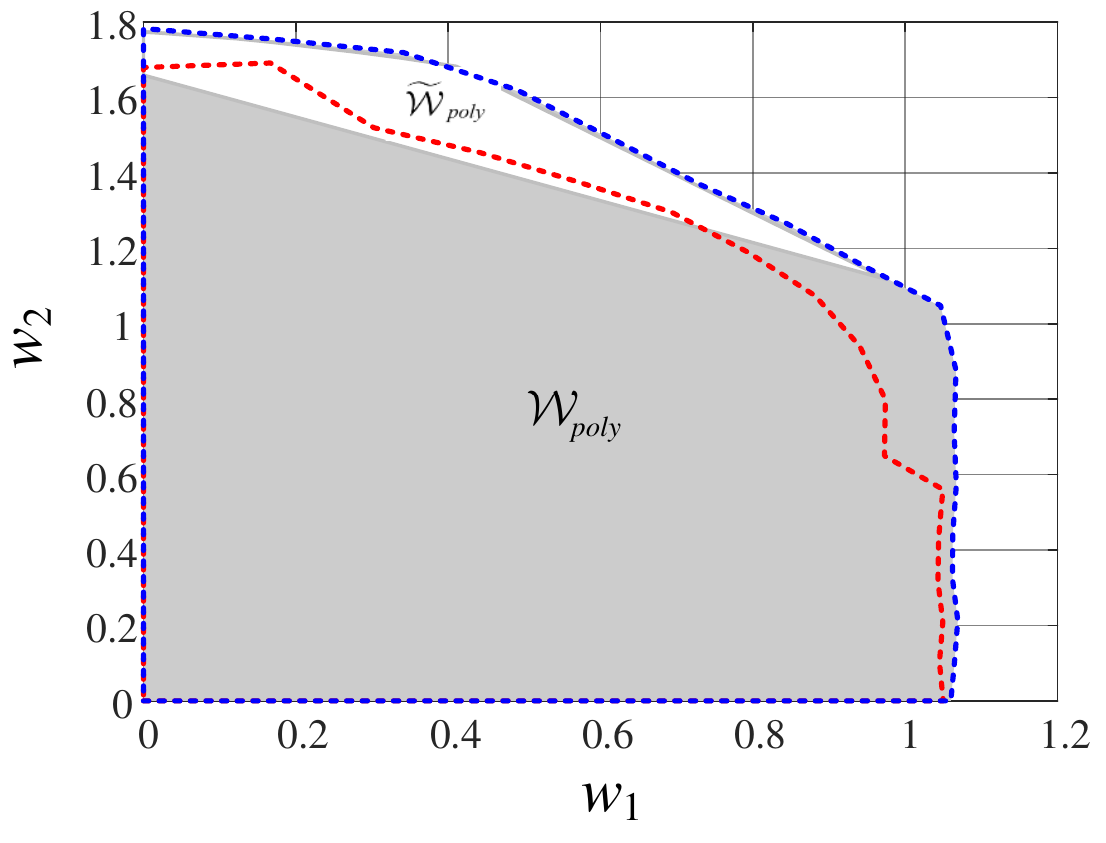}
        \caption{The gray polytope $\mathcal{W}_{poly}$ is an approximation of $\mathcal{W}$ (red dash line). It is obtained by removing the output $\tilde{\mathcal{W}}_{poly}$ of Algorithm \ref{alg:remove-inexact} (white polytope) from the output $\mathcal{W}'_{poly}$ of Algorithm \ref{alg:approximate-socp} (the  outside blue dash line).}
        \label{fig:algorithm2}
\end{figure}

Even though Algorithm 1 can output the accurate SOCP-relaxed dispatchable region, as we can see in Fig. \ref{fig:actual-relax}, there is still a gap between the actual dispatchable region $\mathcal{W}$ and the relaxed one $\mathcal{W}'$. If the renewable generator output $(w_1,w_2)$ lies in the gap area, there is actually no feasible dispatch that satisfies power flow equation \eqref{eq:dist-flow} and safety limits \eqref{eq:pq_limits}-\eqref{eq:vl_limits}. Thus, using the SOCP-relaxed region as a guidance will threaten power system security. In this paper, Algorithm \ref{alg:remove-inexact} is developed to further remove the nondispatchable points. As in Fig. \ref{fig:algorithm2}, the $\tilde{\mathcal{W}}_{poly}$ (white area) generated by Algorithm \ref{alg:remove-inexact} is removed and the resulting region $\mathcal{W}_{poly}$ (grey area) is closer to the actual region (red dash line). This shows the great potential of the proposed algorithm in improving the accuracy of dispatchable region in a distribution system. The operational risk under the obtained region $\mathcal{W}_{poly}$ and the SOCP-relaxed region $\mathcal{W}'$ will be compared later in TABLE \ref{tab:compare}.

\subsection{Impact of different factors} \label{sec:numerical-2}
In the following, we test the impact of two factors (adjustable capability of controllable generators $[\underline{p}_i,\overline{p}_i],\forall i$ and current limit $\overline \ell$) on the shape of the dispatchable region and the performance of the proposed algorithm.

\begin{table}[t]
        \renewcommand{\arraystretch}{1.3}
        \renewcommand{\tabcolsep}{1em}
        \centering
        \caption{Parameters of generators in Cases L and H}
        \label{tab:generator2}
        \begin{tabular}{ccccc}
                \hline 
                Generator &  \multicolumn{2}{c}{Case L} & \multicolumn{2}{c}{Case H} \\ No. & $\underline{p}_i$ (p.u.) & $\overline{p}_i$ (p.u.) & $\underline{p}_i$ (p.u.) & $\overline{p}_i$ (p.u.)\\
                \hline
                G1 & 0.4 & 0.5 & 0 & 0.6\\
                G2 & 0.3 & 0.4 & 0 & 0.4\\
                G3 & 0.4 & 0.5 & 0 & 0.6\\
                G4 & 0.4 & 0.5 & 0 & 0.5\\
                G5 & 0.4 & 0.5 & 0 & 0.6\\
                \hline
        \end{tabular}
\end{table}
\begin{figure}[t]
        \centering
        \includegraphics[width=1.0\columnwidth]{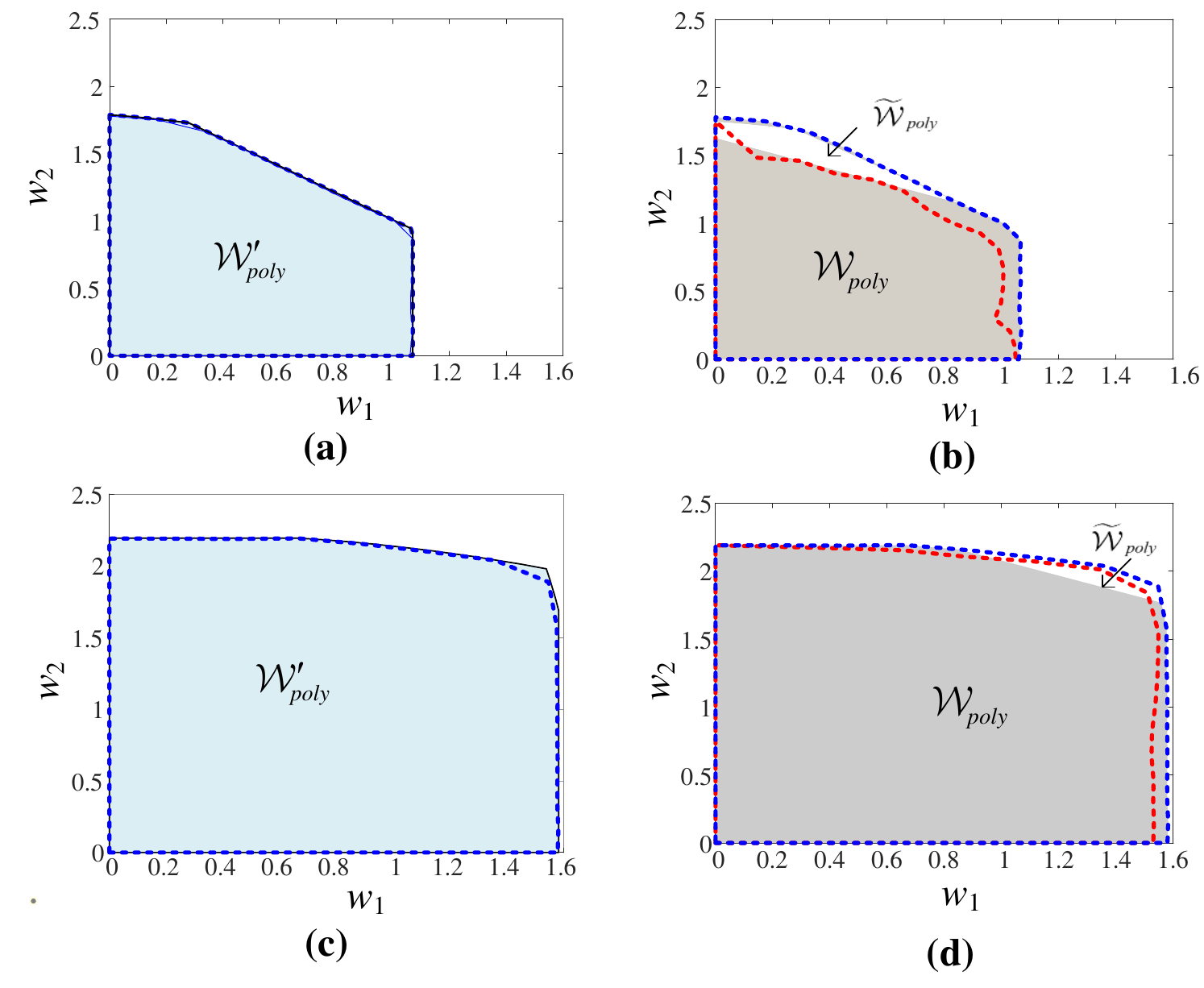}
        \caption{Left: the $\mathcal{W}'_{poly}$ returned by Algorithm \ref{alg:approximate-socp} for Case L (subfigure (a)) and Case H (subfigure (c)). Right: the $\mathcal{W}_{poly}$ returned by Algorithm \ref{alg:remove-inexact} for Case L (subfigure (b)) and Case H (subfigure (d)).}
        \label{fig:comparison}
\end{figure}
\begin{figure}[t]
	\includegraphics[width=0.8\columnwidth]{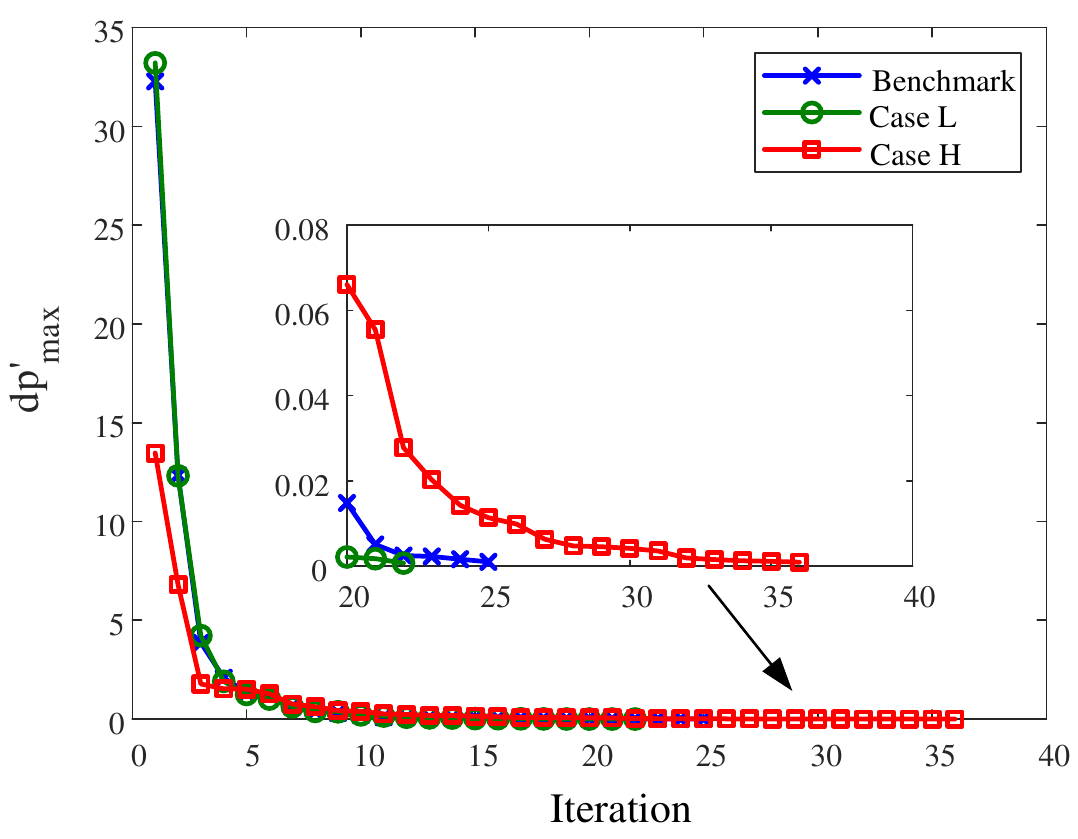}
	\centering
	\caption{The change of $\mathrm{dp}'_{max}$ over iterations of Algorithm \ref{alg:approximate-socp}.}
	\label{fig:dpmax}
\end{figure}
To show how $[\underline{p}_i,\overline{p}_i]$ influences the dispatchable region, we test three cases: (1) \textbf{Benchmark}, which has the same setting as in Section \ref{sec:numerical-1}. (2) \textbf{Case L}, where the generators have less adjustable capability than the benchmark. (3) \textbf{Case H}, where the generators have more adjustable capability than the benchmark. The parameters of the Cases L and H are given in TABLE \ref{tab:generator2}. The regions returned by Algorithm \ref{alg:approximate-socp} ($\mathcal{W}'_{poly}$) and Algorithm \ref{alg:remove-inexact} ($\mathcal{W}_{poly}$) are given in Fig. \ref{fig:comparison}. Subfigures (a), (b) are for Case L and subfigures (c), (d) are for Case H. The changes of $\mathrm{dp}'_{max}$ under three cases are recorded in Fig. \ref{fig:dpmax}.

As shown in Figs. \ref{fig:alg1} and \ref{fig:comparison}, Algorithm \ref{alg:approximate-socp} can always output the accurate SOCP-relaxed dispatchable region, i.e., $\mathcal{W}'_{poly}=\mathcal{W}'$. The final dispatchable regions (grey area) returned by the proposed algorithms are much closer to the actual ones compared with the SOCP-relaxed regions. In addition, as the adjustable capability of generators decreases, the system's ability to accommodate volatile renewable power becomes weaker, and thus, the dispatchable region becomes smaller. We also find that with a weaker adjustable capability, the actual dispatchable region $\mathcal{W}$ is more likely to be nonconvex and to differ more from the SOCP-relaxed region. The difference between the red dash line and the blue dash line in Fig. \ref{fig:comparison}(b) is more significant than that in Fig. \ref{fig:comparison}(d).
In the future power systems, more renewable generators are replacing the controllable generators, so the use of an SOCP-relaxed dispatchable region is not accurate enough.
Therefore, the proposed Algorithms \ref{alg:approximate-socp}-\ref{alg:remove-inexact} to remove the nondispatchable points will be helpful. 

\begin{table}[t]
        \renewcommand{\arraystretch}{1.3}
        \renewcommand{\tabcolsep}{1em}
        \centering
        \caption{Comparison of three cases.}
        \label{tab:compare}
        \begin{tabular}{cccccc}
                \hline 
                 & \textbf{FR}($\mathcal{W}'$) & \textbf{FR}($\mathcal{W}_{poly}$) & Reduction & Time(s) \\
                 \hline
                 Benchmark & 10.4\% & 4.5\% & 56.73\% & 289.84 \\
                 Case L & 15.7\% & 8.7\% & 44.59\% & 291.01\\
                 Case H & 3.5\% & 2.5\% & 28.57\% & 724.60\\
                \hline
        \end{tabular}
\end{table}

In Fig. \ref{fig:dpmax}, the $\mathrm{dp}'_{max}$ under all three cases decrease towards zero when Algorithm \ref{alg:approximate-socp} terminates. The computational times are 289.84s (Benchmark), 291.01s (Case L), and 724.60s (Case H), respectively, showing that our algorithm is efficient. Moreover, we randomly generate 2000 points $(w_1,w_2)$ in the SOCP-relaxed dispatchable region $\mathcal{W}'$ and the final obtained region $\mathcal{W}_{poly}=\mathcal{W}'_{poly}  \backslash \tilde{\mathcal{W}}_{poly}$, and calculate the failure rate defined as
\begin{align}\label{eq:failurerate}
    \textbf{FR}(\mathcal{S})=\frac{\mbox{No. of points}~w \in \mathcal{S}~\mbox{that is nondispatchable}}{\mbox{No. of points}~w\in \mathcal{S}}
\end{align}
The failure rates under three cases are summarized in TABLE \ref{tab:compare}. In all three cases, the proposed method can greatly reduce the failure rate, and the reduction is more than 50\% under benchmark. This can help better ensure system security. Moreover, we can find that in a system with relatively small adjustable capability, the reduction is more significant.

\begin{figure}[t]
	\includegraphics[width=1.0\columnwidth]{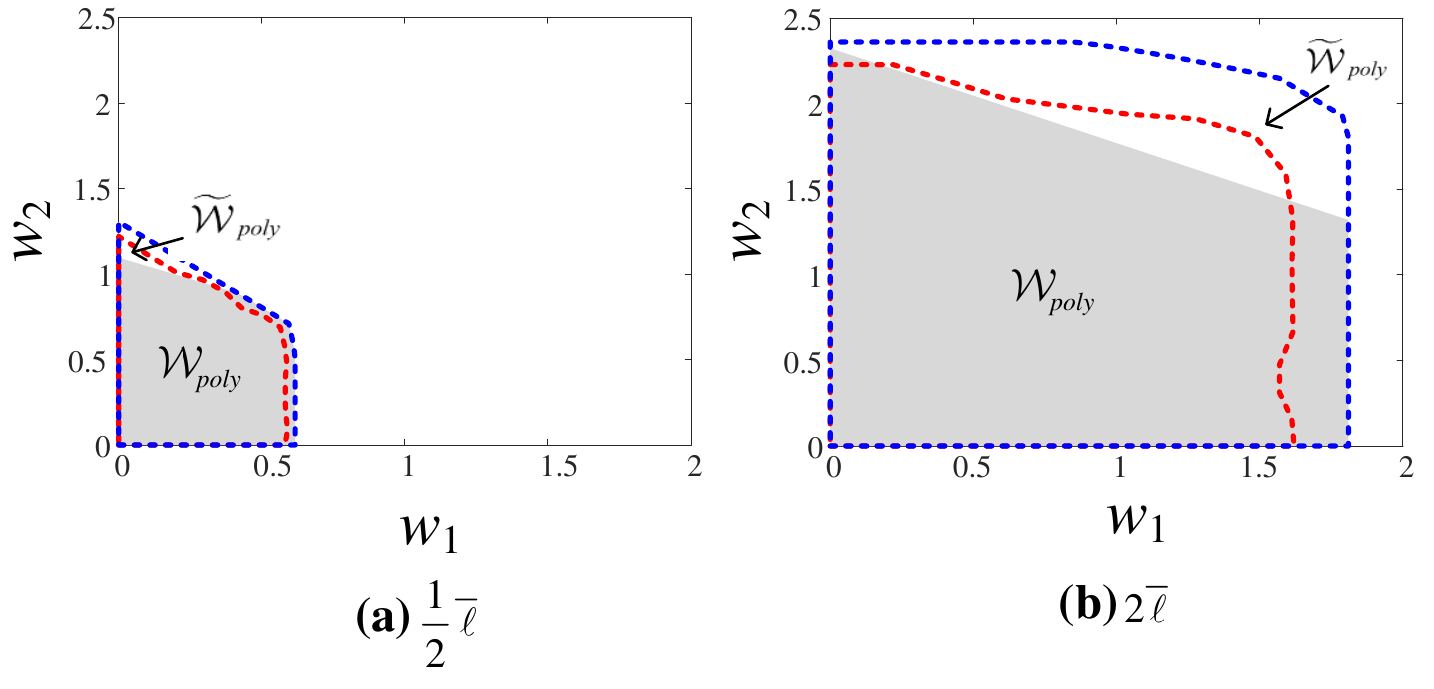}
	\centering
	\caption{The $\mathcal{W}_{poly}$ returned by Algorithm 2 under two current limits.}
	\label{fig:current}
\end{figure}

Furthermore, we test the impact of current limit by running two other cases where we halve and double the $\overline \ell$, respectively. The obtained region are shown in Fig. \ref{fig:current}. The computation times are both less than 500s, which is acceptable. A more stringent line-flow limit results in a smaller dispatchable region and also a greater deviation between the relaxation region $\mathcal{W}'$ and the exact region $\mathcal{W}$.

%We run Algorithm \ref{alg:remove-inexact} multiple times as proposed in Section \ref{sec:methods}. The obtained convex polytopes (whose union serves as an approximation of $\mathcal{\tilde W}$) are removed from the output $\mathcal{W}'_{poly}$ by Algorithm \ref{alg:approximate-socp}, leaving an approximation of the dispatchable region $\mathcal{W}$. Fig. \ref{fig:comparison} compares this approximation with $\mathcal{W}$ obtained by checking sampled points (close to the actual $\mathcal{W}$), which shows the ability of the proposed algorithms to provide a simple and relatively accurate approximation. 

\subsection{Comparison with other methods}
\begin{figure*}[t]
        \centering
        \includegraphics[width=1.8\columnwidth]{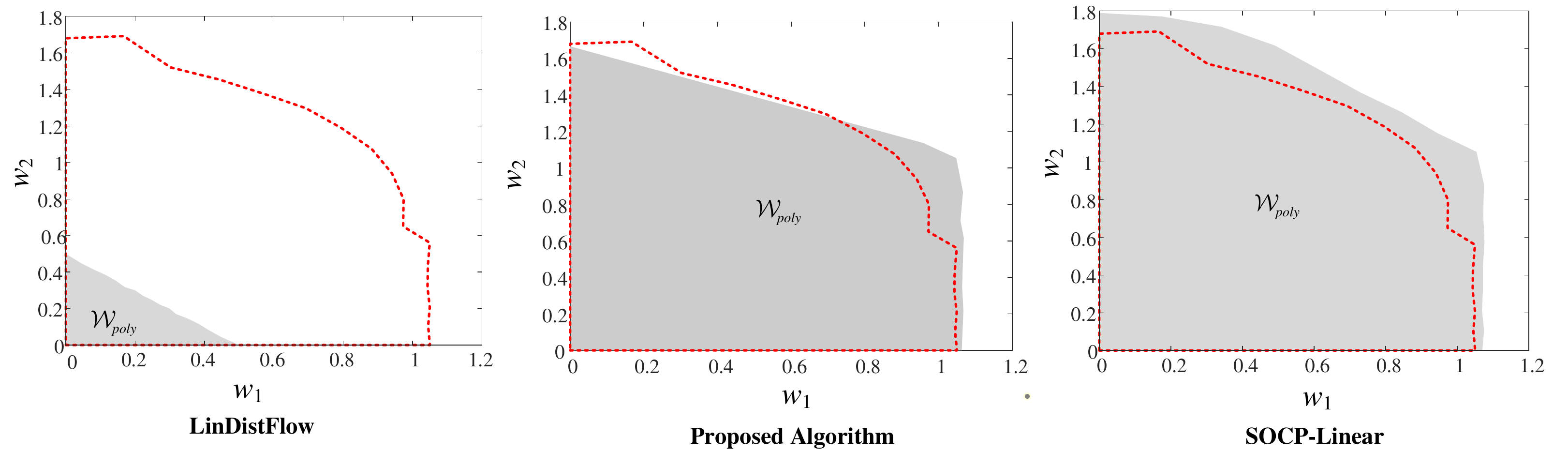}
        \caption{The regions returned by Linearized DistFlow model (left), the proposed algorithms (middle), and the polyhedral approximation of SOCP-relaxed model (right).}
        \label{fig:differentmethod}
\end{figure*}
We then compare the performance of the proposed algorithms with two well-known approaches based on (1) Linearized DistFlow model \cite{baran1989optimal} (denoted as \textbf{LinDistFlow}); (2) Polyhedral approximation of the SOCP-relaxed model \cite{chen2018energy} (denoted as \textbf{SOCP-Linear}). These two models are linear programs so that the adaptive constraint generation algorithm in \cite{wei2015real} can be applied to generate the region. The results are shown in Fig. \ref{fig:differentmethod}. Theoretically, the LinDistFlow region can be an inner/outer approximation or a region that intersects the actual dispatchable region. In the benchmark case, the LinDistFlow 
region is very small and conservative. The SOCP-Linear region is always an outer approximation of the actual region. We can see that it is very close to the SOCP-relaxed region $\mathcal{W}'$. To better illustrate the results under three approaches, we calculate the failure rate \eqref{eq:failurerate} and the missing rate (\textbf{MR}) defined below.
\begin{align}
    \textbf{MR}(\mathcal{S})=\frac{\mbox{No. of points}~ w \in \mathcal{W}~\mbox{but}~\notin \mathcal{S}}{\mbox{No. of points}~w \in \mathcal{W}}
\end{align}
The failure rate and missing rate under three approaches are compared in TABLE \ref{tab:threemethod}. We can find that, in this simulation case, the LinDistFlow region is an inner approximation so its failure rate is zero. However, it has a very high missing rate, meaning that the region is too conservative. The SOCP-Linear region is always an outer approximation so its missing rate is zero, but its failure rate is high. The proposed method can achieve a good balance between ensuring security and reducing conservatism.
\begin{table}[t]
        \renewcommand{\arraystretch}{1.3}
        \renewcommand{\tabcolsep}{1em}
        \centering
        \caption{Comparison of three methods.}
        \label{tab:threemethod}
        \begin{tabular}{cccccc}
                \hline 
                 & \textbf{FR}($\mathcal{W}_{poly}$) & \textbf{MR}($\mathcal{W}_{poly}$) \\
                 \hline
                 LinDistFlow & 0\% & 91.1\% \\
                 Proposed Method & 4.5\% & 2.7\% \\
                 SOCP-Linear & 11.8\% & 0\%\\
                \hline
        \end{tabular}
\end{table}

\section{Conclusion}\label{sec:conclusion}
In this paper, we develop an improved approximation of the renewable generation dispatchable region in radial distribution networks. First, a nonconvex optimization problem is formulated to describe the dispatchable region. The nonconvex problem is then relaxed to a convex SOCP. An SOCP-based projection algorithm (Algorithm 1) is proposed to generate the accurate SOCP-relaxed dispatchable region under certain conditions. In addition, a heuristic method (Algorithm 2) is developed to remove the SOCP-inexact region from the region obtained above. Therefore, the final region can better approximate the actual nonconvex dispatchable region. Our main findings are:
\begin{itemize}
    \item The proposed method can reduce the operational risk (quantified by failure rate) by more than 50\% compared with the SOCP-relaxed region.
    \item The proposed method has a greater potential in the future power system with fewer controllable units and thus weaker adjustable capability.
    \item Compared with existing approaches (LinDistFlow and SOCP-Linear), the proposed method achieves a better tradeoff between security and conservatism.
\end{itemize}

This paper provides an innovative perspective for constructing the dispatchable region: While the existing literature can only generate convex regions, the proposed algorithm can generate nonconvex approximations. For future work, we aim to improve the accuracy of the proposed algorithms by properly setting the initial points for heuristic searching.

%In this paper, we develop an improved polytopic approximation to the dispatchable region of renewable power generations in radial distribution networks. First, an optimization problem to check the AC power-flow solvability and safety limits satisfaction is formulated to describe the dispatchable region. The nonconvex problem is then relaxed to a convex SOCP, whose dual problem (also an SOCP) was proven to attain strong duality. An SOCP based projection algorithm is proposed that can generate the exact SOCP-relaxed dispatchable region under certain conditions. In addition, a heuristic method is developed to remove the SOCP-inexact region from the region obtained above. Therefore, the final region can better approximate the nonconvex actual dispatchable region. We have to admit that the heuristic method has no performance guarantee, but it indeed provides an innovative perspective for constructing the dispatchable region: While the existing literature can only generate convex regions, the proposed algorithm can generate nonconvex approximations. For future work, we plan to improve accuracy of the proposed algorithms by properly setting the initial points for heuristic searching. 

\section*{Appendix. Constant parameters}

This appendix provides in full detail the constant matrices, vectors, and numbers used in Section \ref{sec:methods}. 

\subsection{Equation \eqref{eq:opt-feasibility}: $A_f$, $B_f$, $\gamma_f$, $A_s$, $\gamma_s$}

The vector $x=(p,q, v,\ell, P,Q)$ is arranged in the order explained in Section \ref{sec:model}. Let $C\in \{-1,0,1\}^{(N+1)\times N}$ be the incidence matrix of the radial network, with its element at the $k$-th row, $j$-th column:
\begin{IEEEeqnarray}{rCl}
C_{kj} &=& \begin{cases}
1,\qquad\text{if}~k=i~\text{for line}~i\rightarrow j \\
-1,\quad \text{if}~k=j~\text{for line}~i\rightarrow j \\
0,\qquad \text{otherwise.}
\end{cases}\nonumber
\end{IEEEeqnarray}
Removing the first row of $C$, we get the reduced incidence matrix $\overline C \in \{-1,0,1\}^{N\times N}$. Define diagonal matrices $R:=\text{diag}(r_{ij}, \forall i\rightarrow j)$ and $X:=\text{diag}(x_{ij}, \forall i\rightarrow j)$. Denote the $N \times N$ all-zero matrix as $\mathbf{O}_{N}$, identity matrix as $I_N$, and $N$-dimensional all-zero column vector as $0_N$. We have: 
\begin{IEEEeqnarray}{rCl}
A_f &=& \begin{bmatrix}
I_N &\mathbf{O}_{N}&\mathbf{O}_{N} & -R & -\overline C & \mathbf{O}_{N} \\
\mathbf{O}_N &I_N&\mathbf{O}_{N} &  -X & \mathbf{O}_{N} & -\overline C \\
\mathbf{O}_N &\mathbf{O}_N &\overline C^\intercal & \left(R^2 \!+\!X^2\right) & -2R & -2X
\end{bmatrix}\nonumber\\
\nonumber\\
\gamma_f &=& \left[0_N^\intercal,~
0_N^\intercal,~
v_0,~ 0_{N-1}^{\intercal}\right]^\intercal.\nonumber
\end{IEEEeqnarray}
Moreover, we define: 
\begin{IEEEeqnarray}{rCl}
B'_f &=&\left[
I_N,~\mathbf{O}_{N},~\mathbf{O}_{N} \right]^\intercal
\nonumber
\end{IEEEeqnarray}
and let $B_f$ be a submatrix of $B'_f$ that contains only the columns corresponding to the nodes $i$ with nonzero renewable generation $w_i$. 
Define column vectors $\overline v := (\overline v_i, ~\forall i=1,...,N)$, $\underline v := (\underline v_i, ~\forall i=1,...,N)$, similarly $\overline p$, $\underline p$, $\overline q$, $\underline q$, and $\overline \ell =(\overline \ell_{ij}, ~\forall i\rightarrow j)$. To write inequalities \eqref{eq:pq_limits}\eqref{eq:vl_limits} as $A_s x + \gamma_s \leq 0$, we need: 
\begin{IEEEeqnarray}{rCl}
A_s &=& \begin{bmatrix}
I_N & \mathbf{O}_{N} & \mathbf{O}_{N}& \mathbf{O}_{N} & \mathbf{O}_{N} & \mathbf{O}_{N} \\
-I_N & \mathbf{O}_{N} & \mathbf{O}_{N}&  \mathbf{O}_{N} & \mathbf{O}_{N} & \mathbf{O}_{N}\\
 \mathbf{O}_{N} & I_N &\mathbf{O}_{N}& \mathbf{O}_{N} & \mathbf{O}_{N} & \mathbf{O}_{N} \\
 \mathbf{O}_{N} &-I_N & \mathbf{O}_{N}&  \mathbf{O}_{N} & \mathbf{O}_{N} & \mathbf{O}_{N}\\
  \mathbf{O}_{N} & \mathbf{O}_{N}& I_N &\mathbf{O}_{N} & \mathbf{O}_{N} & \mathbf{O}_{N} \\
 \mathbf{O}_{N} &\mathbf{O}_{N}&-I_N &   \mathbf{O}_{N} & \mathbf{O}_{N} & \mathbf{O}_{N}\\
\mathbf{O}_{N}& \mathbf{O}_{N} & \mathbf{O}_{N} & I_N & \mathbf{O}_{N} & \mathbf{O}_{N} \\
\mathbf{O}_{N}& \mathbf{O}_{N} & \mathbf{O}_{N} & -I_N & \mathbf{O}_{N} & \mathbf{O}_{N}
\end{bmatrix},~ 
\gamma_s =\begin{bmatrix}
-\overline p \\
\underline p \\
-\overline q \\
\underline q \\
-\overline v \\
\underline v \\
-\overline \ell \\
0_N
\end{bmatrix}.\nonumber
\end{IEEEeqnarray}

\subsection{Equation \eqref{eq:opt-feasibility-relaxed}: $A_{y}$, $b_{y}$, $c_{q}$, $\gamma_{q}$}

To make \eqref{eq:opt-feasibility:soc-substitute}--\eqref{eq:opt-feasibility:soc} the same as: 
\begin{IEEEeqnarray}{rCl}
 \left\| \begin{bmatrix}
2P_{ij}  \\
2Q_{ij} \\
v_i - \ell_{ij}
\end{bmatrix} \right\|_2
&\leq& v_i + \ell_{ij} +z_{q,ij}, \quad  \forall i\rightarrow j  \nonumber
\end{IEEEeqnarray}
we need $A_{y}$, $b_{y}$, $c_{q}$, $\gamma_{q}$ as follows:
\begin{itemize}
    \item For all $i \rightarrow j$, $A_{y,ij}$ is $3\times (6N)$ sparse matrix with all elements zero except its element at the first row, $(4N+j)$-th column equal to $2$; at the second row, $(5N+j)$-th column equal to 2; at the third row, $(2N+i)$-th column equal to $1$ (if $i\neq 0$), and $(3N+j)$-th column equal to $-1$.   
    \item For all $i\rightarrow j$ except $0\rightarrow 1$, $b_{y,ij}$ is a three-dimensional column vector of all zeros; $b_{y,01} = \left[0, 0, v_0 \right]^\intercal$.
    
    \item For all $i\rightarrow j$, $c_{q,ij}$ is a $(6N)$-dimensional row vector of all zeros except its $(2N+i)$-th (if $i\neq 0$) and $(3N+j)$-th elements both equal to $1$. 
    \item $\gamma_{q,ij}=0$ for all $i\rightarrow j$ except $0\rightarrow 1$; $\gamma_{q,01} = v_0$. 
\end{itemize}

\ifCLASSOPTIONcaptionsoff
\newpage
\fi
        
        \bibliographystyle{IEEEtran}
        \bibliography{IEEEabrv,mybib}

\end{document}